\documentclass[11pt]{article}
\usepackage[english]{babel}
\usepackage[left=1in,right=1in,top=1in,bottom=1.5in]{geometry} 
\usepackage[skip=5pt, indent=10pt]{parskip}
\usepackage{amssymb, amsmath, amsthm}
\usepackage{hyperref}
\usepackage{xcolor}
\usepackage{nicefrac}
\usepackage[margin=1cm]{caption}
\usepackage{natbib}
\usepackage{bbm}
\usepackage[scr]{rsfso}
\usepackage{algorithm}
\usepackage{algpseudocode}
\usepackage{times}
\usepackage{graphicx}
\usepackage{stackrel}
\usepackage{cases}

\definecolor{darkred}{HTML}{880000}
\definecolor{darkblue}{HTML}{000088}

\hypersetup{
  colorlinks,
  linkcolor={darkred},
  citecolor={darkblue}
}

\newcommand{\myendproof}{\hfill$\square$}

\renewcommand{\leq}{\leqslant}

\renewcommand{\geq}{\geqslant}

\renewcommand{\tilde}{\widetilde}

\newcommand{\eps}{\varepsilon}

\newcommand{\dd}{{\mathrm{d}}}

\newcommand{\sfp}{\mathsf{p}}

\newcommand{\NN}{\mathsf{NN}}

\newcommand{\E}{\mathbb E}
\newcommand{\N}{\mathbb N}
\newcommand{\p}{\mathbb P}
\newcommand{\R}{\mathbb R}
\newcommand{\Z}{\mathbb Z}
\newcommand{\1}{\mathbbm 1}

\newcommand{\sR}{\mathscr R}

\newcommand{\bj}{{\mathbf j}}
\newcommand{\bk}{{\mathbf k}}

\newcommand{\cB}{{\mathcal{B}}}

\newcommand{\cH}{{\mathcal{H}}}
\newcommand{\cK}{{\mathcal{K}}}

\newcommand{\cN}{{\mathcal{N}}}
\newcommand{\cO}{{\mathcal{O}}}
\newcommand{\cP}{{\mathcal{P}}}
\newcommand{\cQ}{{\mathcal{Q}}}

\newcommand{\cS}{{\mathcal{S}}}

\newcommand{\cU}{{\mathcal{U}}}

\newcommand{\floor}[1]{{\lfloor #1 \rfloor}}
\newcommand{\ceil}[1]{\lceil #1 \rceil}

\def\esssup{\operatornamewithlimits{esssup}}

\newcommand{\gelu}{\mathrm{GELU}}
\newcommand{\id}{\mathrm{id}}

\newcommand{\integral}[1]{\int\limits_{#1}}

\newtheorem{Th}{Theorem}[section]
\newtheorem{Lem}[Th]{Lemma}
\newtheorem{Def}[Th]{Definition}

\newtheorem{As}[Th]{Assumption}

\title{Simultaneous Approximation of the Score Function\\ and Its Derivatives by Deep Neural Networks}

\author{
Konstantin Yakovlev\thanks{HSE University, Russian Federation, kdyakovlev@hse.ru}
\and
Nikita Puchkin\thanks{HSE University, Russian Federation, npuchkin@hse.ru}
}

\date{}

\begin{document}

\maketitle

\begin{abstract}
    We present a theory for simultaneous approximation of the score function and its derivatives, enabling the handling of data distributions with low-dimensional structure and unbounded support.
    Our approximation error bounds match those in the literature while relying on assumptions that relax the usual bounded support requirement.
    Crucially, our bounds are free from the curse of dimensionality. 
    Moreover, we establish approximation guarantees for derivatives of any prescribed order, extending beyond the commonly considered first-order setting.
\end{abstract}

\section{Introduction}

Score estimation, the task of learning the gradient of the log density, has become a crucial part of generative diffusion models \citep{song2019generative,song2021scorebased}.
These models achieve state-of-the-art performance in a wide range of domains including images, audio and video synthesis \citep{dhariwal2021diffusion,kong2021diffwave,ho2022video}.
To sample from the desired distribution, one needs to have an accurate score function estimator along the Ornstein-Uhlenbeck process.
In the context of diffusion models the score estimation is done through the minimization of \emph{denoising score matching} loss function over the class of neural networks \citep{song2021scorebased,vincent2011connection,oko2023diffusion}.
Another recipe for score estimation is \emph{implicit score matching} proposed by \cite*{hyvarinen2005estimation}.
The proposed objective includes not only the score function, but also its Jacobian trace.

A crucial research question is to determine the \emph{iteration complexity} of the distribution estimation given inaccurate score function.
The convergence theory of diffusion models has received much attention in the recent years.
Some works \citep{de2022convergence,chen2023sampling,benton2024nearly,li2024adapting} study SDE-based samplers under the assumption that the score estimator is $L^2$-accurate.
However, several recent works \citep{yang2023nearly,li2024accelerating,li2024towards,li2024sharp,huang2025fast,li2025faster} demonstrate that  ODE-based samplers can achieve faster convergence rate, provided that not only the score function but also its Jacobian matrix are well-approximated.
Notably, \cite{huang2025fast} further requires boundedness of the second derivatives of the estimated score.
A similar issue arises in the convergence analysis of Langevin dynamics under the inaccurate drift assumption \citep{huggins2017quantifying,majka2020nonasymptotic}.
In these works, the approximate drift is assumed to be Lipschitz continuous, a condition that is satisfied if the drift has bounded derivatives on the entire space.

Recent works considering denoising score matching \citep{oko2023diffusion,chen2023score,tang2024adaptivity,azangulov2024convergence,yakovlev2025generalization} primarily focus on quantification of the $L^2$-error leaving the study of higher-order derivatives out of the scope.
To our knowledge, only \cite{shen2023differentiable} addresses this issue providing rigorous guarantees for approximation of first-order derivatives.
However, the problem of approximating higher-order score derivatives remains open.
Furthermore, \cite{shen2023differentiable} assumes that the target distribution has a bounded support while generative diffusion models usually deal with data corrupted with Gaussian noise.

Recent advances in analysis of differentiable neural networks have demonstrated their ability to approximate not only functions but also their derivatives.
For example, Sobolev training \citep{czarnecki2017sobolev} and physics-informed neural networks \citep{raissi2018hidden,raissi2019physics,mishra2023estimates} leverage derivative information to improve generalization ability of the resulting estimators.
Furthermore,  \cite{guhring2021approximation,de2021approximation,abdeljawad2022approximations,belomestny2023simultaneous} establish that differentiable neural networks can universally approximate functions and their derivatives in a range of smoothness classes, including Sobolev and Hölder spaces.
Their findings suggest that neural networks are well-suited for tasks requiring higher-order smoothness, such as score estimation through implicit score matching.
The only drawback of these works is that they focus on full-dimensional settings and face the curse of dimensionality. In contrast, we are interested in scenarios where data lies near a low-dimensional manifold, which is typical for real-world data \citep{bengio2013representation,pope2021the,brown2023verifying}.
Although some works have explored approximation capabilities of neural networks with a rectified linear unit (ReLU) activation function under the assumption that the input data lies exactly on a low-dimensional manifold or has a small Minkowski dimension \citep{schmidt2019deep,chen2019efficient,nakada20,jiao2023deep}, the problem of simultaneously approximating functions and their derivatives in a noisy setup, which is common in analysis of generative diffusion models, remains underexplored.

We also note that score estimation along with its high-order derivatives is crucial when estimating the density ridge \citep{genovese2014nonparametric,sasaki2018mode}.
Specifically, to estimate the density ridge, one must have a precise estimate of the score function and its derivatives up to order three in the neighbourhood of the ridge.
Moreover, a common data model for density ridge estimation assumes a hidden manifold structure: the target distribution is generated by first sampling a latent variable supported on an unknown low-dimensional manifold and then adding full-dimensional noise \citep{genovese2014nonparametric}.

This naturally leads to the following research question:

\begin{quote}
    \emph{
    Can a feedforward neural network simultaneously approximate the score function and all its derivatives up to any prescribed order, while operating in the presence of full-dimensional noise and avoiding the curse of dimensionality?}
\end{quote}

\medskip
\noindent
\textbf{Contribution.}\quad
We present approximation theory for simultaneous approximation of the score function and its derivatives up to an arbitrarily large predefined order.
Our theory covers the case of data distributions obtained by first applying a smooth mapping to a uniformly distributed input on a low‑dimensional unit cube and then convolving the result with Gaussian noise.
Specifically, we consider distributions for which a random vector $X_0 \in \R^D$ is generated according to the model:
\begin{align*}
    X_0 = g^*(U) + \sigma Z,
\end{align*}
where $U \sim \mathrm{Un}([0, 1]^d)$ and $Z \sim \cN(0, I_D)$ are independent.
Here, the unknown mapping $g^* : [0, 1]^d \to \R^D$ belongs to a certain smoothness class.
Notably, the effective dimension $d$ may be significantly smaller than the ambient dimension $D$, which is a common feature of real-world data arising in tasks such as density ridge estimation \citep{genovese2014nonparametric} and diffusion-based generative modeling \citep{yakovlev2025generalization}.
A key advantage of our approach is that it does not impose the restriction of bounded support on the target distribution, in contrast to prior work such as \cite{shen2023differentiable}.

Extending the approximation framework originally developed for rectified linear unit (ReLU) neural networks \citep{yakovlev2025generalization}, we consider Gaussian error linear units (GELU) \citep{hendrycks2016gaussian}, an infinitely smooth alternative to ReLU.
We leverage the results from \cite{yakovlev2025gelu}, which showed that feedforward neural networks with GELU activation functions can approximate elementary operations including multiplication, exponential functions, and division.
We also carefully track the dependence on both the ambient dimension $D$ and the noise scale $\sigma$ in the derivation of the architecture of the neural network.
Note that capturing this dependence is crucial in high-dimensional settings.

\medskip
\noindent
\textbf{Paper structure.}\quad
The rest of the paper is organized as follows.
In Section \ref{sec:prelim_not}, we introduce all necessary notation and definitions.
Section \ref{sec:main_results} presents our main result and highlights its distinction from related work.
A full proof of the main result is provided in Section \ref{sec:proof_main_res}, while auxiliary lemmas and technical results are deferred to the appendix.

\medskip
\noindent
\textbf{Notation.}\quad We denote by $\Z_+$ the set of non-negative integers.
The tuple $\bk \in \Z_+^d$ for some $d \in \N$ is referred to as a multi-index, which is always displayed in bold.
We also write $|\bk| = k_1 + k_2 + \ldots + k_d$, $\bk! = k_1! \cdot k_2! \cdot \ldots \cdot k_d!$, and for a vector $v \in \R^d$ we define $v^\bk = v_1^{k_1} v_2^{k_2} \dots v_d^{k_d}$.
For a subset $\Omega \subseteq \R^d$ and a function $f : \Omega \to \R$ we denote by
\begin{align*}
    \partial^\bk f = \frac{\partial^{|\bk|} f}{\partial x_1^{k_1} \partial x_2^{k_2} \dots \partial x_d^{k_d}}
\end{align*}
its weak derivative.
For any $1 \leq i \leq j \leq d$ and any $v \in \R^d$ the expression $v_{i:j}$ stands for $(v_i, v_{i + 1}, \dots, v_j)$.
For any $R > 0$ and a vector $v \in \R^d$, a Euclidean ball of radius $R$ centered at $v$ is denoted as $\cB(v, R)$.
In addition, we write $f \lesssim g$ if $f = \cO(g)$.
If $f \lesssim g$ and $g \lesssim f$, then we write $f \asymp g$.
The notation $\tilde{\cO}(\cdot)$ is defined similar to $\cO(\cdot)$, except that it ignores logarithmic terms.
Throughout the paper, we often replace the expression for $\min\{a, b\}$ and $\max\{a, b\}$ with $a \vee b$ and $a \wedge b$ correspondingly.
For any $x > 0$ we define $\log(x) = \ln(x \vee e)$.

\section{Preliminaries and notations}
\label{sec:prelim_not}

\textbf{Norms.}\quad For a vector $v$ we denote its Euclidean norm by $\|v\|$, the maximal absolute value of its entries by $\|v\|_\infty$, and the number of its non-zero entries by $\|v\|_0$.
Similarly, $\|A\|_\infty$ and $\|A\|_0$ denote the maximal absolute value of entries of $A$ and the number of its non-zero entries, respectively.
Finally, for a function $f : \Omega \to \R^d$ and non-negative weight function $w : \Omega \to \R$, we define
\begin{align*}
    \|f\|_{L^\infty(\Omega)} = \esssup_{x \in \Omega}\|f(x)\|,
    \quad \|f\|_{L^2(\Omega, w)} = \left\{\integral{\Omega}\|f(x)\|^2 w(x) \, \dd x\right\}^{1 / 2} .
\end{align*}
When the set $\Omega$ is clear from context, we often write $\|f\|_{L^2(w)}$ to avoid the abuse of notation.

\medskip
\noindent
\textbf{Smoothness spaces.}\quad
To characterize regularity of functions in our analysis, we introduce the Sobolev space and the Hölder class.
Their formal definitions are given below.

\begin{Def}[Sobolev space]
    Let $r \in \N$ and $\Omega \subseteq \R^r$ be an open set.
    Define the Sobolev space $W^{k, \infty}(\Omega)$ with $k \in \Z_+$ as
    \begin{align*}
        W^{k, \infty}(\Omega) = \{f \in L^\infty(\Omega) : \partial^\bk f \in L^\infty(\Omega) \quad \text{for every } \bk \in \Z_+^r \text{ with } |\bk| \leq k \},
    \end{align*}
    where $L^\infty(\Omega)$ is the Lebesgue space.
    We define the Sobolev seminorm on $W^{k, \infty}(\Omega)$ as
    \begin{align*}
        |f|_{W^{k, \infty}(\Omega)} = \max_{\substack{\bk \in \Z_+^r, \; |\bk| = k}}\|\partial^\bk f\|_{L^\infty(\Omega)}.
    \end{align*}
    Based on the seminorm, introduce the Sobolev norm on $W^{k, \infty}(\Omega)$
    \begin{align*}
        \|f\|_{W^{k, \infty}(\Omega)} = \max_{0 \leq m \leq k}|f|_{W^{m, \infty}(\Omega)} .
    \end{align*}
\end{Def}

\begin{Def}[Hölder class]
    \label{def:holder_class}
    For given $\beta > 0$ we let $\floor{\beta}$ to be the largest integer strictly less than $\beta$.
    Let also the function space $C^\floor{\beta}(\Omega)$ contain the functions $f : \Omega \to \R$ which have bounded and continuous derivatives up to order $\floor{\beta}$, that is,
    \begin{align*}
        C^\floor{\beta} = \left\{f : \Omega \to \R : \max_{\bk \in \Z_+^r, \; |\bk| \leq \floor{\beta}} \|\partial^\bk f\|_{L^\infty(\Omega)} < \infty \right\} .
    \end{align*}
    The Hölder class consists of those functions in $C^\floor{\beta}(\Omega)$ whose derivatives of order $\floor{\beta}$ are $(\beta - \floor{\beta})$-Hölder-continuous.
    Formally, we define
    \begin{align*}
        \cH^\beta(\Omega, \R, H) = \left\{f \in C^\floor{\beta}(\Omega) : \max_{\substack{\bk \in \Z_+^r \\ |\bk| \leq \floor{\beta}}} \|\partial^\bk f\|_{L^\infty(\Omega)} \leq H,
        \; \max_{\substack{\bk \in \Z_+^r \\ |\bk| = \floor{\beta}}} \sup_{x\not= y \in \Omega}\frac{|\partial^\bk f(x) - \partial^\bk f(y)|}{1 \wedge \|x - y\|_\infty^{\beta - \floor{\beta}}} \leq H \right\} .
    \end{align*}
    The Hölder class for vector-valued functions is defined component-wise. 
    That is, a mapping $h : \Omega \to \R^m$, for some $m \in \N$, belongs to $\cH^\beta(\Omega, \R^m, H)$ if and only if each coordinate function $h_i$ lies in $\cH^\beta(\Omega, \R, H)$ for $1 \leq i \leq m$.
\end{Def}

\medskip
\noindent
\textbf{Neural networks.}\quad
In this paper we focus on feed-forward neural networks with the activation function
\begin{align*}
    \gelu(x) = x \cdot \Phi(x), \quad \Phi(x) = \frac{1}{\sqrt{2\pi}}\integral{-\infty}^x e^{-t^2 / 2} \dd t .
\end{align*}
The choice of the activation function is motivated by the fact that it is infinitely smooth, and all its derivatives are uniformly bounded (see Lemma \ref{lem:gelu_seminorms_bound}).
Next, for a vector $b = (b_1, \dots, b_r) \in \R^r$ we define the shifted activation function $\gelu_b : \R^r \to \R^r$ as
\begin{align*}
    \gelu_b(x) = (\gelu(x_1 - b_1), \dots, \gelu(x_r - b_r)), \quad x = (x_1, \dots, x_r) \in \R^r .
\end{align*}
For some $L \in \N$ and a vector $W = (W_0, W_1, \dots, W_L) \in \N^{L + 1}$, a neural network of depth $L$ and architecture $W$ is a function $f : \R^{W_0} \to \R^{W_L}$ such that
\begin{align}
    \label{eq:feed_forward_nn_def}
    f(x) = -b_L + A_L \circ \gelu_{b_{L - 1}} \circ A_{L - 1} \circ \gelu_{b_{L - 2}} \circ \dots \circ A_2 \circ \gelu_{b_1} \circ A_1 \circ x ,
\end{align}
where $A_j \in \R^{W_{j} \times W_{j - 1}}$ is a weight matrix and $b_j \in \R^{W_j}$ is a bias vector for all $j \in \{1, \dots, L\}$.
The maximum number of neurons of each layer is denoted as $\|W\|_\infty$ and is referred to as the width of the neural network.
We define the class of neural networks of the form \eqref{eq:feed_forward_nn_def} with at most $S$ non-zero weights and the weight magnitude $B$ as follows:
\begin{align}
    \label{eq:nn_class_def}
    \NN(L, W, S, B) = \left\{f \text{ of the form } \eqref{eq:feed_forward_nn_def} : \sum_{j = 1}^L (\|A_j\|_0 + \|b_j\|_0) \leq S, \; \max_{1 \leq j \leq L} \|A_j\|_\infty \vee \|b_j\|_\infty \leq B \right\} .
\end{align}

\section{Main results}
\label{sec:main_results}

We begin by formulating the main assumption of our work, which allows us to mitigate the curse of dimensionality.

\begin{As}
    \label{asn:relax_man}
    Let $g^* \in \cH^\beta([0, 1]^d, \R^D, H)$, for some $H > 0$ and $d, D \in \N$, with $\|g^*\|_{L^\infty([0, 1]^d)} \leq 1$.
    Let $\sigma > 0$, and suppose that the data sample $X_0 \in \R^D$ is generated by the model:
    \begin{align*}
        X_0 = g^*(U) + \sigma \, Z,
    \end{align*}
    where $U \sim \mathrm{Un}([0, 1]^d)$ and $Z \sim \cN(0, I_D)$ are independent random elements.
\end{As}

Informally, Assumption \ref{asn:relax_man} models distributions concentrated near a low-dimensional manifold, perturbed by isotropic Gaussian noise.
This is motivated by the common observation that real-world data often has a small intrinsic dimension despite residing in a high-dimensional space \citep{levina2004maximum,pope2021the}.
This assumption is milder than those in related works, as it does not require the image of $g^*$ to have a positive reach \citep{tang2024adaptivity,azangulov2024convergence} or its distribution to have a lower-bounded density with respect to the volume measure \citep{oko2023diffusion,tang2024adaptivity,azangulov2024convergence}.
Similar relaxed assumptions have been successfully used to analyze the convergence rates of generative adversarial networks \citep{schreuder2021statistical,stephanovitch2024wasserstein} and diffusion models \citep{yakovlev2025generalization}, and they play a crucial role in alleviating the curse of dimensionality.

In view of Assumption \ref{asn:relax_man}, the data distribution density $\sfp^*(y)$ has an explicit form
\begin{align}
    \label{eq:p_star_def}
    \sfp^*(y) = (\sqrt{2\pi}\sigma)^{-D} \int_{[0, 1]^d} \exp\left\{-\frac{\|y - g^*(u)\|^2}{2\sigma^2}\right\} \, \dd u, \quad y \in \R^D .
\end{align}
Therefore, score function also can be explicitly written
\begin{align}
    \label{eq:s_star_f_star_def}
    s^*(y)
    = \nabla \log \sfp^*(y)
    = -\frac{y}{\sigma^2} + \frac{f^*(y)}{\sigma^2},
    \quad f^*(y) = \frac{\integral{[0, 1]^d}g^*(u)\exp\left\{-\frac{\|y - g^*(u)\|^2}{2\sigma^2}\right\} \dd u  }{\integral{[0, 1]^d} \exp\left\{-\frac{\|y - g^*(u)\|^2}{2\sigma^2}\right\} \dd u } .
\end{align}
The derived result motivates us to consider GELU neural networks with tailored architecture.
We now present the main results of the paper, focusing on the quantitative expressive power of the neural network class defined in \eqref{eq:nn_class_def} in approximating the true score function.

\begin{Th}[approximation of the true score function]
    \label{thm:approx_main}
    Grant Assumption \ref{asn:relax_man}.
    Also assume that $\eps \in (0, 1)$ is sufficiently small in the sense that it satisfies
    \[
        \eps^\beta
        \leq \frac{\floor{\beta}!}{H d^\floor{\beta} \sqrt{D}}\left(1 \wedge \frac{C_1 \sigma^2}{\sqrt{D}m^2\left(\log(1 / \eps) + \log\left(mD\sigma^{-2}\right) \right)} \right)
    \]
    and
    \[
        (H \vee 1)^2 P(d, \beta)^2 D m^2(\log(1 / \eps) + \log(mD\sigma^{-2}))\eps \leq C_2 \sigma^2 ,
    \]
    where $C_1$ and $C_2$ are absolute positive constants.
    Then for any $m \in \N$ there exists a score function approximation $\bar{s} \in \cS(L, W, S, B)$ of the form
    \begin{align*}
        \bar{s}(y) = -\frac{y}{\sigma^2} + \frac{\bar{f}(y)}{\sigma^2}, \quad y \in \R^D,
    \end{align*}
    which satisfies
    \begin{align*}
        (i) &\quad
        \max_{1 \leq l \leq D}\max_{\bk \in \Z_+^D, \, |\bk| \leq m} \|\partial^\bk[\bar{s}_l - s^*_l] \|^2_{L^2(\sfp^*)}
        \lesssim \sigma^{-4|\bk| - 8} e^{\cO(|\bk|\log |\bk|)} D^2\eps^{2\beta} \log^2(1 / \eps) \log^2 \left(mD\sigma^{-2} \right), \\
        (ii) &\quad \max_{1 \leq l \leq D} |\bar{f}_l|_{W^{k, \infty}(\R^D)} \leq \sigma^{-2k}\exp\left\{\cO \Big( k^2 \log\left(mD\log(1 / \eps)\log\left(\sigma^{-2} \right) \right) \Big) \right\}, \quad \text{for all } 0 \leq k \leq m .
    \end{align*}
    Moreover, $\bar{f}$ has the following configuration:
    \begin{align*}
        &L \lesssim \log\big(mD\sigma^{-2}\log(1 / \eps) \big),
        \quad \log B \lesssim m^{85}D^8 \log^{26}(mD\sigma^{-2})\log^{21}(1 / \eps), \\
        &\|W\|_\infty \vee S
        \lesssim \eps^{-d} D^{16 + P(d, \beta)} m^{132 + 17 P(d, \beta)}\sigma^{-48 - 4P(d, \beta)} \left( \log(mD\sigma^{-2}) \log(1 / \eps) \right)^{38 + 4 P(d, \beta)},
    \end{align*}
    where $P(d, \beta) = \binom{d + \floor{\beta}}{d}$.

\end{Th}

A full proof of Theorem \ref{thm:approx_main} is given in Section \ref{sec:proof_main_res}.
We now compare our result with that of \cite{shen2023differentiable}.
For this purpose, let us set $m = 1$ and assume that the precision parameter $\eps$ is sufficiently small.
Under these conditions, Theorem \ref{thm:approx_main} guarantees the mean square error of $\widetilde\cO(\eps^{2\beta})$ for both the score function and all of its derivatives up to order $m$.
This is achieved by a neural network with $\widetilde\cO(\eps^{-d})$ non-zero parameters, thereby mitigating the curse of dimensionality.
This matches the dependence on $\eps$ obtained in \cite{shen2023differentiable}.

However, there are key differences in the assumptions of the two results.
First, \cite{shen2023differentiable} assumes that the data is supported in a neighbourhood of radius $\rho$ of a sufficiently regular manifold and has bounded support,
while our Assumption \ref{asn:relax_man} allows for unbounded data distributions with low-dimensional structure.
Note that the image of $g^*$ is not necessarily a manifold.
Second, the analysis in \cite{shen2023differentiable} heavily relies on a general approximation result for functions from the space $C^s$ for some $s \in \N$ (see Definition \ref{def:holder_class}).
A key limitation of this approach is that it does not provide approximation guarantees for derivatives of order higher than one.
In contrast, by leveraging specific properties of the score function under Assumption \ref{asn:relax_man}, we circumvent this restriction.
Since the score function is infinitely smooth in our setting (see \citep[Lemma 4.1]{yakovlev2025generalization}), we establish approximation guarantees not only for the score itself but also for its derivatives of arbitrarily high order.
The most important limitation of \citep{shen2023differentiable} is that the neighborhood radius $\rho$ is required to vanish with the approximation error $\eps$.
This forces the data to lie arbitrarily close to the manifold, allowing their analysis to mimic the case of data supported exactly on it.
For this reason, the setup of \citep{shen2023differentiable} is even more restrictive than the ones of \citep{tang2024adaptivity, azangulov2024convergence}.
In contrast,  or model allows the noise scale $\sigma$ in Assumption \ref{asn:relax_man} to be fixed, directly resolving this limitation and handling a wider class of data distributions.
This work thus establishes the first approximation error bounds for this broad class of distributions that avoid the curse of dimensionality.

We also emphasize that the result derived in Theorem \ref{thm:approx_main} cannot be derived from general approximation theory for functions and their derivatives \citep{belomestny2023simultaneous,de2021approximation}.
Specifically, applying the bound on the Sobolev seminorms of $f^*$ (see \eqref{eq:s_star_f_star_def}) from Lemma \ref{lem:f_circ_f_star_acc}, and invoking Corollary 1 from \cite{belomestny2023simultaneous},
yields an approximation of each component of $s^*$ and its derivatives up to order $m$ with accuracy $\widetilde{\cO}(\eps^{2\beta})$ on $[0, 1]^D$ using a neural network with $\widetilde{\cO}((m D \log(1 / \eps))^D)$ non-zero parameters.
Evidently, when the ambient dimension $D \gtrsim \log(1 / \eps)$, this parameter count exceeds that given in Theorem \ref{thm:approx_main}.
Finally, \cite{belomestny2023simultaneous} does not guarantee a finite Sobolev norm outside the unit cube, limiting its applicability to $L^2$-norm error quantification.

Let us elaborate on the intuition behind claim $(ii)$ of Theorem \ref{thm:approx_main}.
Its proof builds upon Lemma \ref{lem:f_circ_f_star_acc}, which establishes that, for each $1 \leq l \leq D$ and $1 \leq k \leq m$, it holds that $|f^*_l|_{W^{k, \infty}} \leq \sigma^{-2k}\exp\{\cO(k \log k)\}$.
Our bound, however, deteriorate faster with $k$ as $\exp\{\cO(k^2 \log m)\}$.
This faster growth with respect to $k$ is an artifact of our proof.
It originates from applying the general Sobolev norm bound for function composition (Lemma \ref{lem:comp_sob_norm}).
The composition rule inherently accounts for all derivative combinations, introducing a combinatorial overhead that leads to the $k^2$ dependence.
Crucially, this discrepancy does not limit the applicability of our main result.
Thanks to its mild requirement of controlling only a small number of derivatives ($m \leq 2$), Theorem \ref{thm:approx_main} allows us to quantify the approximation error in the convergence analysis of diffusion models \citep{yang2023nearly, li2024accelerating, li2024towards, li2024sharp, huang2025fast, li2025faster}.

It is worth noting that Theorem \ref{thm:approx_main} can be extended to the case where one aims to approximate the score function along the Ornstein-Uhlenbeck forward process \citep{oko2023diffusion,yakovlev2025generalization}.
This extension is possible since our construction does not require approximating derivatives with respect to the time variable.
We emphasize that by setting $m = 1$ in Theorem \ref{thm:approx_main}, we recover the approximation result of the true score function similar to that derived in \cite[Theorem 3.3]{yakovlev2025generalization}.
However, in our case, the dependence on $\sigma$ becomes worse both in terms of the number of non-zero parameters and the approximation accuracy.
We argue that this dependence could be made identical to that of \cite{yakovlev2025generalization} by adapting their argument if one does not need to approximate derivatives of the score function.

\section{Proof of Theorem \ref{thm:approx_main}}
\label{sec:proof_main_res}

The proof of our main approximation theorem requires a new technique, as standard approaches are insufficient under Assumption \ref{asn:relax_man}.

A natural starting point is the method of \cite[Theorem 8]{shen2023differentiable}.
This approach combines the Whitney extension theorem \citep{fefferman2006whitney} with dimensionality reduction via linear projections \citep{baraniuk2009random}, which is applicable in the case of a sufficiently regular, low-dimensional manifold.
However, this path encounters two fundamental obstacles in our setting.
First, the image of $g^*$ is not guaranteed to be a well-conditioned manifold and, therefore, the result of \cite{baraniuk2009random} is inapplicable.
Second, even under a manifold assumption, we will inevitably encounter issues with limited differentiability, since the embedding into a lower-dimensional space only approximately preserves distances.

Therefore, we develop a new constructive approach.
Our proof is inspired by the technique of \cite{yakovlev2025generalization} but requires significant extension to operate in Sobolev norms of arbitrary order.
A key and novel difficulty is ensuring good behavior of the approximation at infinity, which was not addressed by prior constructions \citep{oko2023diffusion,tang2024adaptivity,azangulov2024convergence,yakovlev2025generalization}.
This difficulty is resolved by leveraging GELU network approximations of elementary functions from \cite{yakovlev2025gelu}.

We now present the complete proof of Theorem \ref{thm:approx_main}, decomposing it into several steps to enhance readability.

\noindent
\textbf{Step 1: local polynomial approximation.}\quad
Following \cite{yakovlev2025generalization}, we introduce the score function $s^\circ$ induced by a local polynomial approximation of $g^*$.
Formally, for each $\bj = (j_1, \dots, j_d) \in \{1, \dots, N\}^d$ with $N = \ceil{1 / \eps}$ we denote
\begin{align*}
    u_\bj = \frac{\bj}{N},
    \quad \cU_\bj = \left[\frac{j_1 - 1}{N}, \frac{j_1}{N}\right] \times \left[\frac{j_2 - 1}{N}, \frac{j_2}{N}\right]
    \times \dots \times \left[\frac{j_d - 1}{N}, \frac{j_d}{N}\right].
\end{align*}
Then the local polynomial approximation is defined for all $[0, 1]^d$ as
\begin{align}
    \label{eq:g_circ_def}
    g^\circ(u) = \sum_{\bj \in \{1, \dots, N\}^d} g^\circ_\bj(u),
    \quad g^\circ_\bj(u) = \sum_{\bk \in \Z_+^d, \; |\bk| \leq \floor{\beta}} \frac{\partial^\bk g^*(u_\bj) (u - u_\bj)^\bk}{\bk!}\1(u \in \cU_\bj) .
\end{align}
We denote
\begin{align}
    \label{eq:s_circ_f_circ_demo}
    s^\circ(y) = -\frac{y}{\sigma^2} + \frac{f^\circ(y)}{\sigma^2},
    \quad f^\circ(y) = \frac{\integral{[0, 1]^d}g^\circ(u)\exp\left\{-\frac{\|y - g^\circ(u)\|^2}{2\sigma^2}\right\} \dd u}{
    \integral{[0, 1]^d}\exp\left\{-\frac{\|y - g^\circ(u)\|^2}{2\sigma^2}\right\} \dd u }.
\end{align}
It is not straightforward how to generalize the result of \cite{yakovlev2025generalization} (Proposition G.1) to the case when we need to evaluate the score derivatives.
We therefore propose a simpler approach based on the following construction.
For a parameter $R > 0$, which will be determined later, we define 
\begin{align}
    \label{eq:cK_R_def}
    \cK = \left\{y \in \R^D : \min_{u \in [0, 1]^d}\|y - g^*(u)\| \leq R \right\}.
\end{align}
The next result quantifies the closeness of $f^*$ and $f^\circ$.
\begin{Lem}
\label{lem:f_circ_f_star_acc}
    Grant Assumption \ref{asn:relax_man}.
    Let $\cK$ be a compact set as defined in \eqref{eq:cK_R_def}, and consider $f^\circ$ as described in \eqref{eq:s_circ_f_circ_demo}.
    Let also $\|g^* - g^\circ\|_{L^\infty([0, 1]^d)} \leq 1 \wedge \sigma^2 / (R + 3)$.
    Then for any $k \in \Z_+$ it holds that
    \begin{align*}
        (i) &\quad \max_{1 \leq l \leq D} \left( |f^\circ_l|_{W^{k, \infty}(\R^D)} \vee |f^*_l|_{W^{k, \infty}(\R^D)} \right) \leq 4^{k + 1} k! \cdot \sigma^{-2k}, \\
        (ii) &\quad \max_{1 \leq l \leq D} |f^\circ_l - f^*_l|_{W^{k, \infty}(\cK)}
        \leq  \sigma^{-2k}(1 + (R + 3) / \sigma^2 )e^{3k + 3} k! \cdot \|g^\circ - g^*\|_{L^\infty([0, 1]^d)} .
    \end{align*}
\end{Lem}

\noindent
We move the proof of Lemma \ref{lem:f_circ_f_star_acc} to Appendix \ref{sec:lem_f_circ_f_star_acc_proof}.
Now we note from \eqref{eq:s_circ_f_circ_demo} that for any $\bk \in \Z_+^D$ and $1 \leq l \leq D$
\begin{align}
    \label{eq:Dk_s_circ_s_star_L2_aux}
    \notag
    &\|D^{\bk}s^\circ_l - D^{\bk}s^*_l\|^2_{L^2(\sfp^*)}
    = \sigma^{-4} \|D^{\bk}f^\circ_l - D^{\bk}f^*_l\|^2_{L^2(\sfp^*)} \\
    &\quad \leq \sigma^{-4}|f^\circ_l - f^*_l|^2_{W^{|\bk|, \infty}(\cK)} + 2\left(|f^\circ_l|^2_{W^{|\bk|, \infty}(\R^D)} + |f^*_l|^2_{W^{|\bk|, \infty}(\R^D)}\right) \integral{\R^D \setminus \cK} \sfp^*(y) \dd y .
\end{align}
Assuming that
\begin{align}
    \label{eq:g_circ_g_star_assn}
    \|g^\circ - g^*\|_{L^\infty([0, 1]^d)} \leq 1 \wedge \sigma^2 / (R + 3) ,
\end{align}
which will be verified later in the proof, we apply Lemma \ref{lem:f_circ_f_star_acc} and obtain that
\begin{align}
    \label{eq:s_circ_s_star_acc}
    \notag
    &\|D^{\bk}s^\circ_l - D^{\bk}s^*_l\|^2_{L^2(\sfp^*)} \\
    &\quad \leq \sigma^{-4|\bk|}(|\bk|!)^2 \left[ \sigma^{- 4} (1 + (R + 3) / \sigma^2 )^2e^{6|\bk| + 6} \|g^\circ - g^*\|_{L^\infty([0, 1]^d)}^2 
    + 4^{2|\bk| + 3} \integral{\R^D \setminus \cK} \sfp^*(y) \dd y \right].
\end{align}
In the next steps, we derive the final bound and optimize the value of $R$.
Specifically, the remainder of the proof involves building a GELU network that is close to $s^\circ$.

\noindent
\textbf{Step 2: reduction to approximation of $s^\circ$ on a compact set.}\quad
We seek an approximation $s$ to $s^\circ$ of the form
\begin{align*}
    s(y) = -\frac{y}{\sigma^2} + \frac{f(y)}{\sigma^2}, \quad f \in \NN(L, W, S, B).
\end{align*}
We also let
\begin{align}
    \label{eq:K_infty_def}
    \cK_\infty = [-R_\infty, R_\infty]^D, \quad R_\infty = 5/2 \vee \sup_{y \in \cK}\|y\|_\infty,
\end{align}
where $\cK$ is defined in \eqref{eq:cK_R_def}.
Thus, using similar argument to \eqref{eq:Dk_s_circ_s_star_L2_aux} and the observation that $\cK \subseteq \cK_\infty$, we arrive at
\begin{align}
    \label{eq:s_s_circ_acc}
    \notag
    &\|D^{\bk}s_l - D^{\bk}s^\circ_l\|^2_{L^2(\sfp^*)}
    = \sigma^{-4} \|D^{\bk}f_l - D^{\bk}f^\circ_l\|^2_{L^2(\sfp^*)} \\
    &\quad \leq \sigma^{-4}|f_l - f^\circ_l|^2_{W^{|\bk|, \infty}(\cK_\infty)} + 2\left(|f_l|^2_{W^{|\bk|, \infty}(\R^D)} + 8^{|\bk| + 1}(|\bk|!)^2\sigma^{-4|\bk|} \right) \integral{\R^D \setminus \cK} \sfp^*(y) \dd y,
\end{align}
where the last inequality uses Lemma \ref{lem:f_circ_f_star_acc}.
Therefore, the goal is to build an approximation of $f^\circ$ on the set $\cK_\infty$ with $|f_l|_{W^{|\bk|, \infty}(\R^D)}$ close to $|f^\circ_l|_{W^{|\bk|, \infty}(\R^D)}$.
We note that \cite{yakovlev2025generalization} used a non-differentiable clipping operation to guarantee that the approximation is bounded at infinity. 
We avoid this obstacle by using a more intricate technique that uses a differentiable approximation of clipping operation.

\noindent
\textbf{Step 3: composition of simpler functions.}\quad
The following steps are devoted to approximation of $f^\circ$ on a compact $\cK_\infty$ defined in \eqref{eq:K_infty_def}.
We first recall from \eqref{eq:s_circ_f_circ_demo} that the expression for $f^\circ(y)$ is given by
\begin{align*}
    f^\circ(y)
    = \frac{1}{\sigma^2}\frac{\integral{[0, 1]^d}g^\circ(u) \exp\left\{-\frac{\|y - g^\circ(u)\|^2}{ 2\sigma^2} \right\} \dd u }{\integral{[0, 1]^d} \exp\left\{-\frac{\|y - g^\circ(u)\|^2}{ 2\sigma^2} \right\} \dd u} .
\end{align*}
In view of \eqref{eq:g_circ_def}, we deduce that $f^\circ$ takes the following form:
\begin{align}
    \label{eq:f_circ_def}
    f^\circ(y) = 
    \frac{\sum_{\bj \in \{1, \dots, N\}^d} \integral{\cU_\bj}g^\circ_\bj(u) \exp\left\{-\frac{\|y - g^\circ_\bj(u)\|^2}{ 2\sigma^2} \right\} \dd u }
    {\sum_{\bj \in \{1, \dots, N\}^d }\integral{\cU_\bj} \exp\left\{-\frac{\|y - g^\circ_\bj(u)\|^2}{ 2\sigma^2} \right\} \dd u}.
\end{align}
Note that for arbitrary $\bj \in \{1, \dots, N\}^d$, $y \in \R^D$ and $u \in \cU_\bj$ we have that
\begin{align}
    \label{eq:square_und_exp_repres}
    \notag
    \frac{\|y - g^\circ_\bj(u)\|^2}{2\sigma^2}
    &= \frac{\|y - g^\circ_\bj(u_\bj)\|^2}{2\sigma^2} + \frac{\|g^\circ_\bj(u) - g^\circ_\bj(u_\bj)\|^2}{2\sigma^2}
    - \frac{(y - g^\circ_\bj(u_\bj))^\top (g^\circ_\bj(u) - g^\circ_\bj(u_\bj))}{\sigma^2} \\
    &= V_{\bj, 0}(y) + \frac{\|g^\circ_\bj(u) - g^\circ_\bj(u_\bj)\|^2}{2\sigma^2} + \sum_{\bk \in \Z_+^d, \; 1 \leq |\bk| \leq \floor{\beta}} V_{\bj, \bk}(y) \frac{(u - u_\bj)^\bk}{\bk!},
\end{align}
where for each $\bj \in \{1, \dots, N\}^d$ we introduced the functions
\begin{align}
    \label{eq:v_j_0_k_def}
    V_{\bj, 0}(y) = \frac{\|y - g^\circ_\bj(u_\bj)\|^2}{2 \sigma^2},
    \quad V_{\bj, \bk}(y) = -\frac{1}{\sigma^2}(y - g^\circ_\bj(u_\bj))^\top \partial^\bk g^*(u_\bj),
    \quad \bk \in \Z_+^d, \; 1 \leq |\bk| \leq \floor{\beta}.
\end{align}
This observation suggests that the integrals
\begin{align}
    \label{eq:ints_over_uj}
    \integral{\cU_\bj}g^\circ_\bj(u) \exp\left\{-\frac{\|y - g^\circ_\bj(u)\|^2}{ 2\sigma^2} \right\} \dd u,
    \quad  \integral{\cU_\bj} \exp\left\{-\frac{\|y - g^\circ_\bj(u)\|^2}{ 2\sigma^2} \right\} \dd u
\end{align}
are compositions of $V_{\bj, \bk}$ and $V_{\bj, 0}$ with a function taking
\begin{align*}
    1 + |\{\bk \in \Z_+^d : 1 \leq |\bk| \leq \floor{\beta}\}|
    = |\{\bk \in \Z_+^{d + 1} : |\bk| = \floor{\beta}\}|
    = \binom{d + \floor{\beta}}{d}
\end{align*}
arguments as input, which may be significantly smaller than the ambient dimension $D$.
For notation simplicity, we denote $P(d, \beta) = \binom{d + \floor{\beta}}{d}$.

\noindent
\textbf{Step 4: approximation of low-dimensional arguments.}\quad
We note from \eqref{eq:v_j_0_k_def} that $V_{\bj, \bk}(y)$ can be implemented exactly using a single linear layer.
The approximation result for $V_{\bj, 0}(y)$ is formulated in the following Lemma.

\begin{Lem}
    \label{lem:exp_V_j_0_approx}
    Let $\cK_\infty \subseteq \R^D$ be a compact set given by \eqref{eq:K_infty_def} and the function $V_{\bj, 0}$ be as defined in \eqref{eq:v_j_0_k_def}.
    Then for any $m \in \N$ and $\eps' \in (0, 1)$ there exists a GELU network $\tilde{\varphi}_{\bj, 0} \in \NN(L, W, S, B)$ such that
    \begin{align*}
        (i) &\quad \|\tilde{\varphi}_{\bj, 0} - \exp \circ (-V_{\bj, 0})\|_{W^{m, \infty}(\cK_\infty)} \leq \eps', \\
        (ii) &\quad \|\tilde{\varphi}_{\bj, 0}\|_{W^{m, \infty}(\R^D)} \leq \exp\{\cO(m^4 \log(mDR\sigma^{-2}\log(1 / \eps')))\} .
    \end{align*}
    Furthermore, $\tilde{\varphi}_{\bj, 0}$ has the following configuration:
    \begin{align*}
        &L \asymp \log(mDR\sigma^{-2}) + \log\log(1 / \eps'),
        \quad \|W\|_\infty \vee S \lesssim D^2 m^{14}(\log^2(1 / \eps') + \log^2(mDR\sigma^{-2})), \\
        &\log B \lesssim m^{23}( \log^3(1 / \eps') + \log^3(mDR\sigma^{-2})) ,
    \end{align*}
    where the number of layers is the same for all $\bj \in \{1, \dots, N\}^d$.

\end{Lem}
The proof of Lemma \ref{lem:exp_V_j_0_approx} is deferred to Appendix \ref{sec:lem_exp_V_j_0_approx}.

\noindent
\textbf{Step 5: approximation of integrals.}\quad 
Introduce a mapping $\sR_\bj : \R^D \to \R^{P(d, \beta)}$ as follows:
\begin{align}
    \label{eq:r_j_def}
    \sR_\bj(y) = \left((V_{\bj, \bk}(y) \; : \; \bk \in \Z_+^d,  \; 1 \leq |\bk| \leq \floor{\beta} ), \frac{1}{2\sigma^2}\right)^\top.
\end{align}
In addition, we define a mapping $a_\bj : \cU_\bj \to \R^{P(d, \beta)}$ componentwise by
\begin{align}
    \label{eq:a_j_def}
    a_\bj(u) = \left( ((u - u_\bj)^\bk / \bk! \; : \; \bk \in \Z_+^d,  \; 1 \leq |\bk| \leq \floor{\beta} ), \|g^\circ_\bj(u) - g^\circ_\bj(u_\bj)\|^2 \right).
\end{align}
Therefore, from \eqref{eq:square_und_exp_repres} we deduce that for all $y \in \R^D$ and $u \in \cU_\bj$, it holds that
\begin{align*}
    \frac{\|y - g^\circ(u)\|^2}{2\sigma^2} = V_{\bj, 0}(y) + \sR_\bj(y)^\top a_\bj(u).
\end{align*}
Hence, the integrals from \eqref{eq:ints_over_uj} simplify to
\begin{equation}
\begin{split}
    \label{eq:ints_U_j_exp_Psi}
    \integral{\cU_\bj}g^\circ_\bj(u) \exp\left\{-\frac{\|y - g^\circ_\bj(u)\|^2}{ 2\sigma^2} \right\} \dd u
    &= \exp(-V_{\bj, 0}(y))\Psi_\bj[g^\circ_\bj](y), \\
    \integral{\cU_\bj} \exp\left\{-\frac{\|y - g^\circ_\bj(u)\|^2}{ 2\sigma^2} \right\} \dd u
    &= \exp(-V_{\bj, 0}(y)) \Psi_\bj[1](y),
\end{split}
\end{equation}
where for a function $\psi$ that can yield vector or scalar outputs, we define
\begin{align}
    \label{eq:Psi_def}
    \Psi_\bj[\psi](y) = \integral{\cU_\bj} \psi(u)\exp\left\{ -\sR_\bj(y)^\top a_\bj(u) \right\} \dd u,
    \quad y \in \R^D.
\end{align}
For the notation simplicity, we denote $\Psi_\bj[1]$ to emphasize that $\psi$ represents the constant function equal to one.
The following results offers approximation bounds for $\Psi_\bj[\psi]$, where $\psi$ is a bounded scalar function.
To complete the current step, it remains to approximate the integrals outlined in \eqref{eq:ints_U_j_exp_Psi}, which is addressed in the following lemma.
\begin{Lem}
    \label{lem:Upsilon_j_gelu_approx}
    Let the assumptions of Lemma \ref{lem:exp_V_j_0_approx} and let
    \begin{align*}
        \Upsilon_\bj[\psi](y) = \exp(-V_{\bj, 0}(y)) \cdot \Psi_\bj[\psi](y), \quad y \in \R^D,
    \end{align*}
    where $\psi : \cU_\bj \to \R$ with $\|\psi\|_{W^{0, \infty}(\cU_\bj)} \leq 2$.
    Let also $\eps$ be sufficiently small so that
    \begin{align*}
        \|g^\circ - g^*\|_{L^\infty([0, 1]^d)} \leq 1,
        \quad \frac{5 (H \vee 1)^2 P(d, \beta)^2\sqrt{D}(R + 3)\eps}{2\sigma^2} \leq 1 .
    \end{align*}
    Then, for any $\eps' \in (0, 1)$ and $m \in \N$ there exists a GELU network $\tilde{\Upsilon}_\bj \in \NN(L_\Upsilon, W_\Upsilon, S_\Upsilon, B_\Upsilon)$ such that
    \begin{align*}
        (i) &\quad \|\tilde{\Upsilon}_\bj - \Upsilon_\bj\|_{W^{m, \infty}(\cK_\infty)} \leq \eps^d \eps', \\
        (ii) &\quad \|\tilde{\Upsilon}_\bj\|_{W^{m, \infty}(\R^D)} \leq \exp\{\cO(m^5 \log(mD\sigma^{-2}\log(1 / \eps')))\} \eps^d .
    \end{align*}
    Furthermore, the configuration of $\tilde{\Upsilon}_\bj$ is as follows:
    \begin{align*}
        &L_\Upsilon \asymp \log(mDR\sigma^{-2}) + \log\log(1 / \eps') ,
        \quad \log B_\Upsilon \lesssim m^{35}\log^3(1 / \eps')\log^7(mDR\sigma^{-2}), \\
        &\|W_\Upsilon\|_\infty \vee S_\Upsilon \lesssim D^2 m^{21 + 7P(d, \beta)}(\log(1 / \eps') + \log(mDR\sigma^{-2}))^{3 + P(d, \beta)} .
    \end{align*}
    where the number of layers does not depend on $\bj \in \{1, \dots, N\}^d$ and the function $\psi$.

\end{Lem}
\noindent
The proof of Lemma \ref{lem:Upsilon_j_gelu_approx} is deferred to Appendix \ref{sec:lem:Upsilon_j_gelu_approx_proof}.
Since $\sup_{u \in [0, 1]^d}\|g_\bj(u)\| \leq 2$, then Lemma \ref{lem:Upsilon_j_gelu_approx} implies that there exist GELU neural networks $\{P_{\bj, l}\}_{l = 1}^D$ and the network $Q_\bj$ that belong to the class $\NN(L_\Upsilon, W_\Upsilon, S_\Upsilon, B_\Upsilon)$ satisfying
\begin{align*}
    \|Q_\bj - \Upsilon_\bj[1]\|_{W^{m, \infty}(\cK_\infty)} \vee \max_{1 \leq i \leq D} \|P_{\bj, l} - \Upsilon_\bj[g^\circ_{\bj, l}]\|_{W^{m, \infty}(\cK_\infty)} \leq \eps^d\eps'.
\end{align*}
Therefore, by Lemma \ref{lem:paral_nn}, the parallel stacking across $\bj \in \{1, \dots, N\}^d$ implies the existence of GELU networks $\{\cP_l\}_{l=1}^D$ and $\cQ$, such that
\begin{align}
    \label{eq:p_q_networks_acc_K}
    \left\|\cQ - \sum_{\bj \in \{1, \dots, N\}^d} \Upsilon_\bj[1]\right\|_{W^{m, \infty}(\cK_\infty)}
    \vee \max_{1 \leq i \leq D}\left\|\cP_l - \sum_{\bj \in \{1, \dots, N\}^d} \Upsilon_\bj[g^\circ_{\bj, l}]\right\|_{W^{m, \infty}(\cK_\infty)}
    \leq \eps' .
\end{align}
Furthermore, Lemmata \ref{lem:Upsilon_j_gelu_approx} and \ref{lem:paral_nn} also suggest that
\begin{align}
    \label{eq:Q_P_l_bound_R}
    \|\cQ\|_{W^{m, \infty}(\R^D)} \vee \max_{1 \leq l \leq D}\|\cP_l\|_{W^{m, \infty}(\R^D)}
    \leq \exp\{\cO(m^5 \log(mDR\sigma^{-2}\log(1 / \eps')))\} .
\end{align}
Finally, we conclude from Lemma \ref{lem:paral_nn} that $\cQ, \cP_1, \dots, \cP_D \in \NN(\breve{L}, \breve{W}, \breve{S}, \breve{B})$ with
\begin{align}
    \label{eq:cQ_cP_cfg}
    \notag
    &\breve{L} \asymp \log(mDR\sigma^{-2}) + \log\log(1 / \eps') ,
    \quad \log \breve{B} \lesssim m^{35}\log^3(1 / \eps')\log^7(mDR\sigma^{-2}), \\
    &\|\breve{W}\|_\infty \vee \breve{S} \lesssim \eps^{-d} D^2 m^{21 + 7P(d, \beta)}(\log(1 / \eps') + \log(mDR\sigma^{-2}))^{3 + P(d, \beta)} .
\end{align}

\noindent
\textbf{Step 5: deriving the final approximation.}\quad
The current step aims to perform division of $\cP_l$ and $\cQ$ for each $1 \leq l \leq D$, as stated in the following lemma.

\begin{Lem}
    \label{lem:div_num_denum_gelu_approx}
    Let the compact set $\cK_\infty \subseteq \R^D$ be defined as in \eqref{eq:K_infty_def}, and $f^\circ$ be as given in \eqref{eq:f_circ_def}.
    Let also the GELU networks $\cQ$ and $\{\cP_l\}_{l=1}^D$ be as defined above in \eqref{eq:p_q_networks_acc_K}.
    Then, choosing
    \begin{align*}
        \log(1 / \eps') \asymp m^6 (R^2 + 1)\sigma^{-2} \log(mDR\sigma^{-2})\log(1 / \eps)
    \end{align*}
    guarantees that there exists a GELU network $\bar{f} \in \NN(\bar{L}, \bar{W}, \bar{S}, \bar{B})$ such that $\bar{f} : \R^D \to \R^D$ and
    \begin{align*}
        (i) &\quad \max_{1 \leq l \leq D} \|\bar{f}_l - f^\circ_l\|_{W^{m, \infty}(\cK_\infty)} \leq \eps^\beta , \\
        (ii) &\quad \max_{1 \leq l \leq D} |\bar{f}_l|_{W^{k, \infty}(\R^D)} \leq \sigma^{-2k} \exp\{\cO( k^2\log(mD\log(1 / \eps)\log(R\sigma^{-2})) )\}, \quad \text{for all } 0 \leq k \leq m .
    \end{align*}
    Furthermore, $\bar{f}$ has
    \begin{align*}
        &\bar{L} \lesssim \log(mDR\sigma^{-2}\log(1 / \eps)),
        \quad \log \bar{B} \lesssim m^{53}R^{16}\sigma^{-16}\log^{10}(mDR\sigma^{-2})\log^4(1 / \eps), \\
        &\|\bar{W}\|_\infty \vee \bar{S}
        \lesssim \eps^{-d} D^4 m^{84 + 13 P(d, \beta)} (R\sigma^{-2})^{24 + 2 P(d, \beta)} (\log(1 / \eps) \log(mDR\sigma^{-2}))^{7 + P(d,\beta)} .
    \end{align*}

\end{Lem}

\noindent
We move the proof of Lemma \ref{lem:div_num_denum_gelu_approx} to Appendix \ref{sec:lem_div_num_denum_gelu_approx_proof}.
The next step is to derive the approximation accuracy of the final approximation given by $\bar{s}(y) = -\sigma^{-2}y + \sigma^{-2}\bar{f}(y)$, where $\bar{f}$ is the neural network from Lemma \ref{lem:div_num_denum_gelu_approx}.
For any $\bk \in \Z_+^D$ with $|\bk| \leq m$ and $1 \leq l \leq D$ we obtain from \eqref{eq:s_circ_s_star_acc}, \eqref{eq:s_s_circ_acc}, Lemma \ref{lem:g_circ_L_inf_acc} and the triangle inequality that
\begin{align*}
    &\|\partial^\bk[\bar{s}_l - s^*_l]\|^2_{L^2(\sfp^*)} \\
    &\quad \lesssim \sigma^{-4|\bk|}(|\bk|!)^2 \sigma^{-4}D\eps^{2\beta}(1 + (R + 3) \sigma^{-2})^2 e^{6|\bk|}
    + (\|\bar{f}\|^2_{W^{m, \infty}(\R^D)} + 8^m(m!)^2\sigma^{-4m})\integral{\R^D \setminus \cK}\sfp^*(y) \dd y .
\end{align*}
Applying Lemma \ref{lem:div_num_denum_gelu_approx} we arrive at
\begin{align}
    \label{eq:D_k_s_bar_s_star_aux}
    &\|\partial^\bk[\bar{s}_l - s^*_l]\|^2_{L^2(\sfp^*)} \\
    \notag
    &\quad \lesssim \sigma^{-4|\bk| - 8}\exp\{\cO(|\bk|\log |\bk|)\}(1 + R^2) D\eps^{2\beta} + \exp\{\cO(m^2\log(mDR\sigma^{-2}\log(1 / \eps)))\} \integral{\R^D \setminus \cK} \sfp^*(y) \dd y .
\end{align}
The next Lemma provides the tail probability upper bound for the data distribution.
\begin{Lem}
    \label{lem:p_star_tail}
    Let $\cK$ be the compact set as defined in \eqref{eq:cK_R_def} with $R \geq \sigma\sqrt{D}$.
    Then it holds that
    \begin{align*}
        \integral{\R^D \setminus \cK} \sfp^*(y) \dd y
        \leq \exp\left\{-\frac{1}{16}\left(-\frac{R^2 - D\sigma^2}{D\sigma^2} \wedge \frac{\sqrt{R^2 - D\sigma^2}}{\sigma}\right)\right\} .
    \end{align*}
\end{Lem}

\noindent
The full proof of Lemma \ref{lem:p_star_tail} is provided in Appendix \ref{sec:lem_p_star_tail_proof}.
Hence, using Lemma \ref{lem:p_star_tail} we deduce that taking
\begin{align}
    \label{eq:R_def}
    R = \sigma\sqrt{D} + 16\sigma\left( \sqrt{D\log(\eps_0^{-2\beta} / D)} \vee \log(\eps_0^{-2\beta} / D) \right)
\end{align}
ensures that
\begin{align*}
    \integral{\R^D \setminus \cK} \sfp^*(y) \, \dd y \leq \eps_0^{2\beta} ,
\end{align*}
where $\eps_0$ is defined as follows:
\begin{align}
    \label{eq:eps_zero_def}
    \log(1 / \eps_0) \asymp m^2\log(1 / \eps) + m^2\log(mD\sigma^{-2}) .
\end{align}
Therefore, from \eqref{eq:D_k_s_bar_s_star_aux} and the above definitions given in \eqref{eq:R_def} and \eqref{eq:eps_zero_def} we deduce that
\begin{align*}
    \|\partial^\bk[\bar{s}_l - s^*_l]\|^2_{L^2(\sfp^*)}
    &\lesssim \sigma^{-4|\bk| - 8}\exp\{\cO(|\bk|\log |\bk|)\} D^2\eps^{2\beta} \log^2(1 / \eps) \log^2(mD\sigma^{-2}) \\
    &\quad + \exp\{\cO(m^2\log(mD\sigma^{-2}\log(1 / \eps)))\} \eps_0^{2\beta} \\
    &\lesssim \sigma^{-4|\bk| - 8}\exp\{\cO(|\bk|\log |\bk|)\} D^2\eps^{2\beta} \log^2(1 / \eps) \log^2(mD\sigma^{-2}) .
\end{align*}
Hence, claim $(i)$ holds true.
From \eqref{eq:R_def} and \eqref{eq:eps_zero_def} we deduce that
\begin{align*}
    R \lesssim \sqrt{D} m^2 (\log(1 / \eps) + \log(mD\sigma^{-2})),
\end{align*}
which together with Lemmata \ref{lem:div_num_denum_gelu_approx} implies that $\bar{f}$ belongs to the class $\NN(L, W, S, B)$ with
\begin{align*}
    &L \lesssim \log(mD\sigma^{-2}\log(1 / \eps)),
    \quad \log B \lesssim m^{85}D^8 \log^{26}(mD\sigma^{-2})\log^{21}(1 / \eps), \\
    &\|W\|_\infty \vee S
    \lesssim \eps^{-d} D^{16 + P(d, \beta)} m^{132 + 17 P(d, \beta)}\sigma^{-48 - 4P(d, \beta)} \left( \log(mD\sigma^{-2}) \log(1 / \eps) \right)^{38 + 4 P(d, \beta)} .
\end{align*}
Note that, in contrast to \cite{yakovlev2025generalization}, the value of $R$ here is not proportional to $\sigma$.
This trade-off is necessary to ensure accurate approximation of the score function and its derivatives.
If only the score function itself needs to be approximated, then using the value of $R$ from \cite{yakovlev2025generalization} leads to a milder dependence on $\sigma$ in the resulting configuration.
Recall from \eqref{eq:g_circ_g_star_assn} and Lemma \ref{lem:Upsilon_j_gelu_approx} that we have to make sure that $\eps$ is sufficiently small so that
\begin{align*}
    \|g^\circ - g^*\|_{L^\infty([0, 1]^d)} \leq 1 \wedge \sigma^2 / (R + 3),
    \quad 5 (H \vee 1)^2 P(d, \beta)^2\sqrt{D}(R + 3)\eps \leq 2\sigma^2 .
\end{align*}
Therefore, Lemma \ref{lem:g_circ_L_inf_acc} suggests that it suffices to take $\eps$ to satisfy
\begin{align*}
    &\eps^\beta
    \leq \frac{\floor{\beta}!}{H d^\floor{\beta} \sqrt{D}}\left(1 \wedge \frac{C_1 \sigma^2}{\sqrt{D}m^2(\log(1 / \eps) + \log(mD\sigma^{-2}))} \right), \\
    &(H \vee 1)^2 P(d, \beta)^2 D m^2(\log(1 / \eps) + \log(mD\sigma^{-2}))\eps \leq C_2 \sigma^2 ,
\end{align*}
where $C_1$ and $C_2$ are some absolute positive constants.
In addition, from \eqref{eq:R_def}, \eqref{eq:eps_zero_def}, and Lemma \ref{lem:div_num_denum_gelu_approx} we find that
\begin{align*}
    \max_{1 \leq l \leq D}|\bar{f}_l|_{W^{k, \infty}(\R^D)}
    \leq \sigma^{-2k} \exp\{\cO( k^2\log(mD\log(1 / \eps)\log(\sigma^{-2})) )\}, \quad \text{for all } 0 \leq k \leq m ,
\end{align*}
which verifies claim $(ii)$.
The proof is now finished.

\bibliographystyle{abbrvnat}
\bibliography{references}

\appendix

\section{Proofs of the auxiliary results}
\label{sec:thm_approx_main_proof}

\subsection{Proof of Lemma \ref{lem:f_circ_f_star_acc}}
\label{sec:lem_f_circ_f_star_acc_proof}
\textbf{Step 1: proving statement $(i)$.}\quad
First, fix an arbitrary $\bk \in \Z_+^D$ and show by induction in $|\bk|$ that for any $1 \leq l \leq D$,
\begin{align}
    \label{eq:Dk_f_star_expr}
    \partial^\bk f_l^*(y) = \frac{\integral{[0, 1]^d}P_{\bk, l}(u_1, \dots, u_{|\bk| + 1}) \exp\left\{ \sum_{j=1}^{|\bk| + 1}\left(\frac{y^\top g^*(u_j)}{\sigma^2} - \frac{\|g^*(u_j)\|^2}{2\sigma^2}\right) \right\} \dd u_1 \dots \dd u_{|\bk| + 1}}
    {\sigma^{2|\bk|} \integral{[0, 1]^d} \exp\left\{ \sum_{j=1}^{|\bk| + 1}\left(\frac{y^\top g^*(u_j)}{\sigma^2} - \frac{\|g^*(u_j)\|^2}{2\sigma^2}\right) \right\} \dd u_1 \dots \dd u_{|\bk| + 1} },
\end{align}
where $P_{\bk, l}$ is some scalar function.
Obviously, the base case when $|\bk| = 0$ holds true for $P_{\mathbf{0}, l}(u_1) = g^*_l(u_1)$.
Assume that the statement is true for all $|\bk| \leq k$ for some $k \in \Z_+$.
Now we perform a transition step.
Note that for all $\bk \in \Z_+^D$ satisfying $|\bk| = k + 1$ one can find $\tilde{\bk}$ with $|\tilde{\bk}| = k$ and $\bk = \tilde{\bk} + \mathbf{e}$,
where $\mathbf{e}$ is a unit vector in a standard basis.
One can show that
\begin{align*}
    D^{\tilde{\bk} + \mathbf{e}} f_l^*(y) = \frac{\integral{[0, 1]^d}P_{\tilde{\bk} + \mathbf{e}, l}(u_1, \dots, u_{k + 2}) \exp\left\{ \sum_{j=1}^{k + 2}\left(\frac{y^\top g^*(u_j)}{\sigma^2} - \frac{\|g^*(u_j)\|^2}{2\sigma^2}\right) \right\} \dd u_1 \dots \dd u_{k + 2}}
    {\sigma^{2(k + 1)}\integral{[0, 1]^d} \exp\left\{ \sum_{j=1}^{k + 2}\left(\frac{y^\top g^*(u_j)}{\sigma^2} - \frac{\|g^*(u_j)\|^2}{2\sigma^2}\right) \right\} \dd u_1 \dots \dd u_{k + 2} },
\end{align*}
where
\begin{align}
    \label{eq:P_k_plus_e_def}
    P_{\tilde{\bk} + \mathbf{e}, l}(u_1, \dots, u_{k + 2}) = P_{\tilde{\bk}}(u_1, \dots, u_{k + 1}) \left( \sum_{j = 1}^{k + 1} \mathbf{e}^\top g^*(u_j) - (k + 1) \cdot \mathbf{e}^\top g^*(u_{k + 2}) \right) .
\end{align}
Therefore, we obtain that
\begin{align*}
    \max_{u_{1: k + 2} \in [0, 1]^d}|P_{\tilde{\bk} + \mathbf{e}, l}(u_{1: k + 2})|
    \leq 2 (k + 1) \max_{u_{1: k + 1} \in [0, 1]^d}|P_{\tilde{\bk}, l}(u_{1: k + 1})| \max_{1 \leq l \leq D} \|g^*_l\|_{W^{0, \infty}([0, 1]^d)}.
\end{align*}
By unrolling the recursion, we deduce that for any $\bk \in \Z_+^D$ it holds that
\begin{align}
    \label{eq:P_k_l_bound}
    \max_{1 \leq l \leq D}\max_{u_{1: |\bk| + 1} \in [0, 1]^d} |P_{\bk, l}(u_{1: |\bk| + 1})|
    \leq 2^{|\bk|} |\bk|! \max_{1 \leq l \leq D}\|g^*_l\|_{W^{0, \infty}([0, 1]^d)}^{|\bk| + 1} .
\end{align}
Hence, from \eqref{eq:Dk_f_star_expr} and Assumption \ref{asn:relax_man} we find that
\begin{align*}
    \max_{1 \leq l \leq D}|f^*_l|_{W^{k, \infty}(\R^D)} \leq 2^k k! \cdot \sigma^{-2k}, \quad \text{for all } k \in \Z_+ .
\end{align*}
Similarly, for $f^\circ$ due to the fact that $\|g^* - g^\circ\|_{L^\infty([0, 1]^d)} \leq 1$, we conclude that
\begin{align*}
    \max_{1 \leq l \leq D}|f^\circ_l|_{W^{k, \infty}(\R^D)} \leq 4^{k + 1} k! \cdot \sigma^{-2k},
    \quad \text{for all } k \in \Z_+ .
\end{align*}
Therefore, statement $(i)$ holds true.

\noindent
\textbf{Step 2: proving statement $(ii)$.}\quad
For each $1 \leq l \leq D$ and $\bk \in \Z_+^D$ let $P_{\bk, l}[g^*](u_{1:|\bk| + 1})$ be the scalar function from \eqref{eq:Dk_f_star_expr} and, similarly, $P_{\bk, l}[g^\circ](u_{1:|\bk| + 1})$.
Let also
\begin{align*}
    Q_k[g](y, u_{1: k + 1}) = \exp\left\{\sum_{j = 1}^{k + 1}\left(\frac{y^\top g(u_j)}{\sigma^2} - \frac{\|g(u_j)\|^2}{2\sigma^2}\right)\right\}, \quad g : [0, 1]^d \to \R^D .
\end{align*}
Then, for any $y \in \cK$ and $\bk \in \Z_+^D$ with $|\bk| \leq m$ we have due to the triangle inequality that
\begin{align}
    \label{eq:f_circ_f_star_abc_dec}
    \sigma^{2|\bk|}|\partial^\bk f^\circ(y) - \partial^\bk f^*(y)|
    \leq (A) + (B) + (C),
\end{align}
where we introduced the terms
\begin{align*}
    (A) &= \frac{\integral{[0, 1]^d}|P_{\bk, l}[g^\circ](u_{1:|\bk| + 1}) - P_{\bk, l}[g^*](u_{1:|\bk| + 1})|Q_{|\bk|}[g^\circ](y, u_{1:|\bk| + 1}) \dd u_{1:|\bk| + 1}}
    {\integral{[0, 1]^d}Q_{|\bk|}[g^\circ](y, u_{1: |\bk| + 1}) \dd u_{1: |\bk| + 1}} , \\
    (B) &= \frac{\integral{[0, 1]^d}|P_{\bk, l}[g^*](u_{1: |\bk| + 1})| \cdot |Q_{|\bk|}[g^\circ](y, u_{1:|\bk| + 1}) - Q_{|\bk|}[g^*](y, u_{|\bk| + 1})| \dd u_{1: |\bk| + 1}}
    {\integral{[0, 1]^d} Q_{|\bk|}[g^\circ](y, u_{1: |\bk| + 1}) \dd u_{1: |\bk| + 1}}, \\
    (C) &= \max_{u_{1: |\bk| + 1} \in [0, 1]^d}|P_{\bk, l}[g^*](u_{1: |\bk| + 1})| \cdot \left|1 - \frac{\integral{[0, 1]^d}Q_{|\bk|}[g^*](y, u_{1: |\bk| + 1}) \dd u_{1: |\bk| + 1} }{\integral{[0, 1]^d} Q_{|\bk|}[g^\circ](y, u_{1: |\bk| + 1}) \dd u_{1: |\bk| + 1}}\right| .
\end{align*}
Now we evaluate each term individually.

\noindent
\textbf{Step 3: bounding term $(A)$.}\quad
We note from \eqref{eq:P_k_plus_e_def} that 
\begin{align*}
    &|P_{\bk + \mathbf{e}, l}[g^*](u_{1: |\bk| + 2}) - P_{\bk + \mathbf{e}, l}[g^\circ](u_{1: |\bk| + 2})| 
    \leq 2(|\bk| + 1)|P_{\bk, l}[g^\circ]| \cdot \|g^\circ - g^*\|_{L^\infty([0, 1]^d)} \\
    &\quad + 2(|\bk| + 1)\|g^*\|_{L^\infty([0, 1]^d)} |P_{\bk, l}[g^*](u_{1: |\bk| + 1}) - P_{\bk, l}[g^\circ](u_{1 : |\bk| + 1})| .
\end{align*}
Using the assumptions of the Lemma and the bound \eqref{eq:P_k_l_bound}, which is also true for $g^\circ$, we deduce that
\begin{align*}
    |P_{\bk + \mathbf{e}, l}[g^*](u_{1: |\bk| + 2}) - P_{\bk + \mathbf{e}, l}[g^\circ](u_{1: |\bk| + 2})|
    &\leq 2(|\bk| + 1) |P_{\bk, l}[g^*](u_{1: |\bk| + 1}) - P_{\bk, l}[g^\circ](u_{1 : |\bk| + 1})| \\
    &\quad + 4(|\bk| + 1)4^{|\bk| + 1} |\bk|! \cdot \|g^\circ - g^*\|_{L^\infty([0, 1]^d)} .
\end{align*}
Thus, unrolling the recursion and taking into account that $|P_{\mathbf{0}, l}[g^*](u_1) - P_{\mathbf{0}, l}[g^\circ](u_1)| \leq \|g^\circ - g^*\|_{L^\infty([0, 1]^d)}$, we arrive at
\begin{align*}
    \max_{1 \leq l \leq D} |P_{\bk, l}[g^*](u_{1: |\bk| + 1}) - P_{\bk, l}[g^\circ](u_{1: |\bk| + 1})|
    \leq 4^{|\bk| + 2}|\bk|! \cdot \|g^\circ - g^*\|_{L^\infty([0, 1]^d)}, \quad \text{for all } \bk \in \Z_+^D .
\end{align*}
Therefore, we have that
\begin{align}
    \label{eq:f_circ_f_star_A_bound}
    (A) \leq 4^{|\bk| + 2} |\bk|! \cdot \|g^\circ - g^*\|_{L^\infty([0, 1]^d)} ,
    \quad \text{for all } y \in \R^D .
\end{align}

\noindent
\textbf{Step 4: bounding terms $(B)$ and $(C)$.}\quad
We note that, for all $y \in \R^D$, the following holds:
\begin{align*}
    &|Q_{|\bk|}[g^*](y, u_{1: |\bk| + 1}) - Q_{|\bk|}[g^\circ](y, u_{1: |\bk| + 1})| \\
    &\quad = Q_{|\bk|}[g^\circ](y, u_{1: |\bk| + 1})
    \left|1 - \exp\left\{\sum_{j = 1}^{|\bk| + 1}\frac{y^\top(g^*(u_j) - g^\circ(u_j))}{\sigma^2} - \frac{\|g^*(u_j)\|^2 - \|g^\circ(u_j)\|^2}{2\sigma^2} \right\}\right|.
\end{align*}
Therefore, the observation that $|\|g^*(u_j)\|^2 - \|g^\circ(u_j)\|^2| \leq 3\|g^\circ - g^*\|_{L^{\infty}([0, 1]^d)}$ and $\sup_{y \in \cK}\|y\| \leq R + 1$ together with the mean value theorem implies that for any $y \in \cK$
\begin{align*}
    &|Q_{|\bk|}[g^*](y, u_{1: |\bk| + 1}) - Q_{|\bk|}[g^\circ](y, u_{1: |\bk| + 1})|
    \\&
    \leq Q_{|\bk|}[g^\circ](y, u_{1: |\bk| + 1}) \frac{(|\bk| + 1)(R + 3)\|g^* - g^\circ\|_{L^\infty([0, 1]^d)}}{\sigma^2}
    \exp\left\{ \frac{(|\bk| + 1)(R + 3)\|g^* - g^\circ\|_{L^\infty([0, 1]^d)}}{\sigma^2} \right\} .
\end{align*}
Using the assumptions of the lemma, we have that for any $y \in \cK$
\begin{align}
    \label{eq:Q_k_star_min_Q_k_circ}
    \frac{|Q_{|\bk|}[g^*](y, u_{1: |\bk| + 1}) - Q_{|\bk|}[g^\circ](y, u_{1: |\bk| + 1})|}{Q_{|\bk|}[g^\circ](y, u_{1: |\bk| + 1})}
    \leq \frac{e^{2|\bk| + 2} (R + 3)\|g^* - g^\circ\|_{L^\infty([0, 1]^d)}}{\sigma^2} .
\end{align}
The obtained bound together with \eqref{eq:P_k_l_bound} implies that for any $y \in \cK$
\begin{align}
    \label{eq:f_circ_f_star_B_bound}
    (B) \leq \frac{ e^{3|\bk| + 2} |\bk|!  (R + 3)\|g^* - g^\circ\|_{L^\infty([0, 1]^d)}}{\sigma^2},
    \quad \text{for all } y \in \cK .
\end{align}
Furthermore, the bounds from \eqref{eq:P_k_l_bound} and \eqref{eq:Q_k_star_min_Q_k_circ} immediately imply that
\begin{align}
    \label{eq:f_circ_f_star_C_bound}
    (C) \leq \frac{e^{3|\bk| + 2} |\bk|! (R + 3)\|g^* - g^\circ\|_{L^\infty([0, 1]^d)}}{\sigma^2},
    \quad \text{for all } y \in \cK .
\end{align}

\noindent
\textbf{Step 5: combining the bounds for $(A)$, $(B)$, and $(C)$.}
\quad
From \eqref{eq:f_circ_f_star_abc_dec}, \eqref{eq:f_circ_f_star_A_bound}, \eqref{eq:f_circ_f_star_B_bound}, and \eqref{eq:f_circ_f_star_C_bound} we deduce that
\begin{align*}
    \max_{1 \leq l \leq D} |f^\circ_l - f^*_l|_{W^{k, \infty}(\cK)}
    \leq  \sigma^{-2k}(1 + (R + 3) / \sigma^2 )e^{3k + 3} k! \cdot \|g^\circ - g^*\|_{L^\infty([0, 1]^d)},
    \quad \text{for all } k \in \Z_+ .
\end{align*}
The proof is finished.

\myendproof

\subsection{Proof of Lemma \ref{lem:exp_V_j_0_approx}}
\label{sec:lem_exp_V_j_0_approx}

\noindent
\textbf{Step 1: approximation of $V_{\bj, 0}$.}\quad
Let $\varphi_{sq}$ be the square operation approximation obtained from Lemma \ref{lem:square_approx} applied with the accuracy parameter $\eps_{sq} \in (0, 1)$ and the smoothness parameter $m \in \N$.
The observation that $\|y\|_\infty \leq 5/2 \vee (R + 1)$ for any $y \in \cK_\infty$ suggests that choosing
\begin{align*}
    \log(1 / \eps_{sq}) \asymp \log(1 / \eps_0) + \log(D R \sigma^{-2})
\end{align*}
ensures that
\begin{align*}
    \max_{1 \leq l \leq D}\left\|\varphi_{sq}((\id)_l - (g^\circ_\bj(u_\bj))_l) -  \| (\id)_l - (g^\circ_\bj(u_\bj))_l \|^2 \right\|_{W^{m, \infty}(\cK_\infty)}
    \leq \exp\{\cO(\log R)\}\eps_{sq}
    \leq \frac{\sigma^2 \eps_0}{D}.
\end{align*}
Here the precision parameters $\eps_0$ and $\eps_{sq}$ will be determined later in the proof.
Therefore, using parallelization argument outlined in Lemma \ref{lem:paral_nn} applied for the networks $\varphi_{sq}(\cdot - (g^\circ_\bj(u_\bj))_l)$
for each $1 \leq l \leq D$ and dividing the last layer parameters on $2\sigma^2$ we obtain that there exists a GELU network $\breve{V}_{\bj, 0}(y)$ satisfying
\begin{align}
    \label{eq:tilde_v_0_acc}
    \notag
    \|\breve{V}_{\bj, 0} - V_{\bj, 0}\|_{W^{m, \infty}(\cK_\infty)}
    &\leq \frac{D}{2\sigma^2} \max_{1 \leq i \leq D}\left\|\varphi_{sq}((\id)_i - (g^\circ_\bj(u_\bj))_i) -  \| (\id)_i - (g^\circ_\bj(u_\bj))_i \|^2 \right\|_{W^{m, \infty}(\cK_\infty)} \\
    &\leq \eps_0 / 2.
\end{align}
Furthermore, from Lemmata \ref{lem:square_approx} and \ref{lem:paral_nn} we deduce that the configuration of $\breve{V}_{\bj, 0}$ is
\begin{align}
    \label{eq:V_breve_j_0_cfg}
    L(\breve{V}_{\bj, 0}) = 2,
    \quad \|W\|_\infty(\breve{V}_{\bj, 0}) \vee S(\breve{V}_{\bj, 0}) \lesssim D,
    \quad \log B(\breve{V}_{\bj, 0}) \lesssim \log(1 / \eps_0) + \log(m D R \sigma^{-2}).
\end{align}
Let $\varphi_{clip}$ be a clipping operation from Lemma \ref{lem:clip_gelu_approx} with the accuracy parameter $\eps_{clip}$, the clipping parameter $5 / 2 \vee (R + 1)$, and the smoothness parameter $m$.
Let also $\varphi_{0, clip}$ we a parallel stacking of $D$ copies of $\varphi_{clip}$ as suggested by Lemma \ref{lem:paral_nn} that approximates a component-wise clipping.
Formally, we have
\begin{align*}
    \max_{1 \leq l \leq D}\|(\varphi_{0, clip})_l - (\id)_l\|_{W^{m, \infty}(\cK_\infty)}
    \leq \eps_{clip}.
\end{align*}
Then, due to Lemma \ref{lem:comp_sob_norm} we have for $\cK_0 = \cK_\infty \cup \varphi_{0, clip}(\R^D)$
\begin{align}
    \label{eq:tilde_v_clip_acc_aux}
    \notag
    &\|\breve{V}_{\bj, 0} \circ \varphi_{0, clip} - \breve{V}_{\bj,0}\|_{W^{m, \infty}(\cK_\infty)} \\
    &\quad \leq \exp\{\cO(m\log(m D))\} \| \breve{V}_{\bj, 0} \|_{W^{m + 1, \infty}(\cK_0)} \eps_{clip} \max_{1 \leq i \leq D}(1 \vee (\|(\id)_i\|_{W^{m, \infty}(\cK_\infty)} + \eps_{clip})^{2m} ).
\end{align}
Lemma \ref{lem:clip_gelu_approx} implies that
\begin{align}
    \label{eq:phi_0_clip_whole}
    \max_{1 \leq l \leq D} \| (\varphi_{0, clip})_l \|_{W^{m, \infty}(\R^D)}
    \leq \exp\{ \cO(m \log m + m\log\log(1 / \eps_{clip}) + \log R  ) \}.
\end{align}
From the definition of $V_{\bj, 0}$ given in \eqref{eq:v_j_0_k_def} and the observation that $\cK_0 \subseteq 2 \cK_\infty$ due to Lemma \ref{lem:clip_gelu_approx}, we obtain that
\begin{align}
    \label{eq:V_j_0_K_0_bound}
    \|V_{\bj, 0}\|_{W^{m, \infty}(\cK_0)}
    \leq \frac{1}{\sigma^2} \| \|\cdot - g^\circ_\bj(u_\bj)\|^2 \|_{W^{m, \infty}(2\cK_\infty)}
    \lesssim \frac{D (R \vee 1)^2}{\sigma^2} .
\end{align}
In addition, we note that one can ensure that \eqref{eq:tilde_v_0_acc} holds for $\cK_0$ with the configuration given in \eqref{eq:V_breve_j_0_cfg}.
Therefore, the triangle inequality yields
\begin{align}
    \label{eq:breve_V_j_0_K0}
    \|\breve{V}_{\bj, 0}\|_{W^{m + 1, \infty}(\cK_0)}
    \leq \|V_{\bj, 0}\|_{W^{m + 1, \infty}(\cK_0)} + \|V_{\bj, 0} - \breve{V}_{\bj, 0}\|_{W^{m + 1, \infty}(\cK_0)}
    \leq \frac{D (R \vee 1)^2}{\sigma^2} .
\end{align}
As a result, \eqref{eq:tilde_v_clip_acc_aux} is evaluated as
\begin{align*}
    \|\breve{V}_{\bj, 0} \circ \varphi_{0, clip} - \breve{V}_{\bj,0}\|_{W^{m, \infty}(\cK_\infty)}
    \leq \exp\{\cO(m\log(m D R \sigma^{-2}) )\} \eps_{clip}.
\end{align*}
Thus, setting
\begin{align}
    \label{eq:eps_clip_v_0}
    \log(1 / \eps_{clip}) \asymp \log(1 / \eps_0) + m \log(mDR\sigma^{-2}),
\end{align}
we ensure that
\begin{align*}
    \|\breve{V}_{\bj, 0} \circ \varphi_{0, clip} - \breve{V}_{\bj,0}\|_{W^{m, \infty}(\cK_\infty)}
    \leq \eps_0 / 2.
\end{align*}
This and \eqref{eq:tilde_v_0_acc} suggest that for $\tilde{V}_{\bj, 0} = \breve{V}_{\bj, 0} \circ \varphi_{0, clip}$ we have
\begin{align}
    \label{eq:tilde_V_acc}
    \|\tilde{V}_{\bj, 0} - V_{\bj, 0}\|_{W^{m, \infty}(\cK_\infty)} \leq \eps_0.
\end{align}
Furthermore, Lemma \ref{lem:comp_sob_norm} together with \eqref{eq:phi_0_clip_whole}, \eqref{eq:breve_V_j_0_K0} and \eqref{eq:eps_clip_v_0} imply that
\begin{align}
    \label{eq:tilde_V_j_0_bound}
    \notag
    \|\tilde{V}_{\bj, 0}\|_{W^{m, \infty}(\R^D)}
    &\leq \exp\{\cO(m\log(m D))\} \|\breve{V}_{\bj, 0}\|_{W^{m, \infty}(\cK_0)} (1 \vee \max_{1 \leq l \leq D}\| (\varphi_{0, clip})_l \|^m_{W^{m, \infty}(\R^D)}) \\
    &\leq \exp\{\cO(m^2\log(mDR\sigma^{-2}\log(1 / \eps_0)))\} .
\end{align}
From \eqref{eq:eps_clip_v_0} and Lemma \ref{lem:clip_gelu_approx} we also find that $\varphi_{0, clip} \in \NN(L_{clip}, W_{clip}, S_{clip}, B_{clip})$ with
\begin{align*}
    L_{clip} = 2,
    \quad \|W_{clip}\|_\infty \vee S_{clip} \lesssim D, 
    \quad \log B_{clip} \lesssim \log(1 / \eps_0) + m \log(mDR\sigma^{-2}) .
\end{align*}
Thus, Lemma \ref{lem:concat_nn} together with \eqref{eq:V_breve_j_0_cfg} suggests that $\tilde{V}_{\bj, 0} \in \NN(\tilde{L}, \tilde{W}, \tilde{S}, \tilde{B})$ with
\begin{align}
    \label{eq:tilde_V_j_0_cfg}
    \tilde{L} &\lesssim 1,
    \quad \|\tilde{W}\|_\infty \vee \tilde{S} \lesssim D^2,
    \quad \log \tilde{B} \lesssim \log(1 / \eps_0) + m \log(mDR\sigma^{-2}) .
\end{align}

\noindent
\textbf{Step 2: composition with the exponential function.}\quad
Now let $\varphi_{exp}$ be a GELU network from Lemma \ref{lem:min_exp_gelu_approx} with the accuracy parameter $\eps_0$, the smoothness parameter $m$, and the parameter $A = 1$.
The triangle inequality suggests that
\begin{align*}
    &\|\varphi_{exp} \circ \tilde{V}_{\bj, 0} - \exp \circ (-V_{\bj, 0}) \|_{W^{m, \infty}(\cK_\infty)} \\
    &\quad \leq \|\varphi_{exp} \circ \tilde{V}_{\bj, 0} - \exp \circ (- \tilde{V}_{\bj, 0}) \|_{W^{m, \infty}(\cK_\infty)} + \| \exp \circ (-\tilde{V}_{\bj, 0}) - \exp \circ (-V_{\bj, 0}) \|_{W^{m, \infty}(\cK_\infty)}.
\end{align*}
From \eqref{eq:V_j_0_K_0_bound}, \eqref{eq:tilde_V_acc} and Lemma \ref{lem:comp_sob_norm} it follows that
\begin{align*}
        \|\varphi_{exp} \circ \tilde{V}_{\bj, 0} - \exp \circ (- \tilde{V}_{\bj, 0}) \|_{W^{m, \infty}(\cK_\infty)}
        \leq \exp\{\cO(m\log(mDR\sigma^{-2})) \} \eps_0.
\end{align*}
and also
\begin{align*}
    \| \exp \circ (-\tilde{V}_{\bj, 0}) - \exp \circ (-V_{\bj, 0}) \|_{W^{m, \infty}(\cK_\infty)}
    \leq \exp\{\cO(m\log(mDR\sigma^{-2}) ) \}\eps_0 .
\end{align*}
Therefore, we obtain that
\begin{align}
    \label{eq:phi_exp_tilde_V_acc}
    \|\varphi_{exp} \circ \tilde{V}_{\bj, 0} - \exp \circ (- V_{\bj, 0}) \|_{W^{m, \infty}(\cK_\infty)}
    \leq \exp\{\cO(m\log(mDR\sigma^{-2}) ) \}\eps_0 .
\end{align}
We also deduce from \eqref{eq:tilde_V_j_0_bound}, Lemma \ref{lem:comp_sob_norm} and Lemma \ref{lem:min_exp_gelu_approx} that
\begin{align}
    \label{eq:phi_exp_tilde_V_bound_R}
    \notag
    \|\varphi_{exp} \circ \tilde{V}_{\bj, 0}\|_{W^{m, \infty}(\R^D)}
    &\leq \exp\{\cO(m\log(m D))\} \| \varphi_{exp} \|_{W^{m, \infty}(\R)} (1 \vee \|\tilde{V}_{\bj, 0}\|_{W^{m,\infty}(\R)}^m) \\
    &\leq \exp\{\cO(m^3\log(mDR\sigma^{-2}\log(1 / \eps_0)))\} .
\end{align}
We now specify the configuration of $\varphi_{exp} \in \NN(L_{exp}, W_{exp}, S_{exp}, B_{exp})$ as suggested by Lemma \ref{lem:min_exp_gelu_approx}:
\begin{align*}
    L_{exp} \lesssim \log m + \log\log(1 / \eps_0),
    \quad \|W_{exp}\|_\infty \vee S_{exp} \lesssim m^6 \log^2(1 / \eps_0),
    \quad \log B_{exp} \lesssim m^{11} \log^3(1 / \eps_0).
\end{align*}
Therefore, Lemma \ref{lem:concat_nn} combined with \eqref{eq:tilde_V_j_0_cfg} implies that $\tilde{\varphi}_{\bj, 0} = \varphi_{exp} \circ \tilde{V}_{\bj, 0} \in \NN(\tilde{L}_\bj, \tilde{W}_\bj, \tilde{S}_\bj, \tilde{B}_\bj)$ with
\begin{align}
    \label{eq:phi_exp_tilde_V_cfg}
    \notag
    &\tilde{L}_\bj \lesssim \log m + \log\log(1 / \eps_0),
    \quad \|\tilde{W}_\bj\|_\infty \vee \tilde{S}_\bj \lesssim D^2 m^6 \log^2(1 / \eps_0), \\
    &\log \tilde{B}_\bj \lesssim m^{11} \log^3(1 / \eps_0) + m \log (mDR\sigma^{-2}).
\end{align}

\noindent
\textbf{Step 3: increasing the number of layers.}\quad
In order to make the number of layers independent of $\bj \in \{1, \dots, N\}^d$, we add identity layers on top of each neural network.
To this end, let $\varphi_{id, \bj}$ be the identity approximation from Lemma \ref{lem:id_deep_gelu_approx} with the precision parameter $\eps_0$,
the number of layers
\begin{align*}
    \log(m + \log(1 / \eps_0)) \vee \max_{\bj' \in \{1, \dots, N\}^d} \tilde{L}_{\bj'} - \tilde{L}_\bj + 1,
\end{align*}
and the scale parameter $K = \|\varphi_{exp} \circ \tilde{V}_{\bj, 0}\|_{W^{m, \infty}(\cK_\infty)}$ .
From \eqref{eq:phi_exp_tilde_V_cfg} we find that the resulting number of layers of each GELU network attains the upper bound for $\breve{L}_\bj$.
Therefore, \eqref{eq:phi_exp_tilde_V_acc}, \eqref{eq:phi_exp_tilde_V_bound_R} and Lemma \ref{lem:comp_sob_norm} imply that for $\tilde{\varphi}_{\bj, 0} = \varphi_{id, \bj} \circ \varphi_{exp} \circ \tilde{V}_{\bj, 0}$ the following holds:
\begin{align*}
    &\|\tilde{\varphi}_{\bj, 0} - \varphi_{exp} \circ \tilde{V}_{\bj, 0} \|_{W^{m, \infty}(\cK_\infty)} \\
    &\quad \leq \exp\{\cO(m \log (mD))\} \| \varphi_{id, \bj} - \id \|_{W^{m, \infty}((\varphi_{exp} \circ \tilde{V}_{\bj, 0}) (\cK_\infty))} (1 \vee \|\varphi_{exp} \circ \tilde{V}_{\bj, 0}\|_{W^{m, \infty}(\cK_\infty)}^m) \\
    &\quad \leq \exp\{\cO(m^4 \log(mDR\sigma^{-2}\log(1 / \eps_0)))\}\eps_0 .
\end{align*}
From this observation and \eqref{eq:phi_exp_tilde_V_acc} it follows that
\begin{align}
    \label{eq:tilde_phi_exp_acc}
    \|\tilde{\varphi}_{\bj, 0} - \exp \circ (-V_{\bj, 0}) \|_{W^{m, \infty}(\cK_\infty)}
    \leq \exp\{\cO(m^4 \log(mDR\sigma^{-2}\log(1 / \eps_0)))\}\eps_0 .
\end{align}
Thus, choosing
\begin{align*}
    \log(1 / \eps_0) \asymp m^4\log(1 / \eps') + m^4\log(mDR\sigma^{-2})
\end{align*}
guarantees that
\begin{align*}
    \|\tilde{\varphi}_{\bj, 0} - \exp \circ (-V_{\bj, 0}) \|_{W^{m, \infty}(\cK_\infty)}
    \leq \eps' .
\end{align*}
In addition, using \eqref{eq:phi_exp_tilde_V_bound_R} in conjunction with Lemmata \ref{lem:comp_sob_norm} and \ref{lem:id_deep_gelu_approx}, we arrive at
\begin{align}
    \label{eq:tilde_phi_exp_zero_R}
    \notag
    \|\tilde{\varphi}_{\bj, 0}\|_{W^{m, \infty}(\R^D)}
    &\leq \exp\{\cO(m\log (mD))\} \|\varphi_{id, \bj}\|_{W^{m, \infty}(\R)} (1 \vee \|\varphi_{exp} \circ \tilde{V}_{\bj, 0}\|_{W^{m, \infty}(\R^D)}^m) \\
    &\leq \exp\{\cO(m^4 \log(mDR\sigma^{-2}\log(1 / \eps')))\}.
\end{align}
Finally, from \eqref{eq:phi_exp_tilde_V_cfg} and Lemma \ref{lem:id_deep_gelu_approx} we deduce that $\tilde{\varphi}_{\bj, 0} \in \NN(L, W, S, B)$ with
\begin{align*}
    &L \asymp \log(mDR\sigma^{-2}) + \log\log(1 / \eps'),
    \quad \|W\|_\infty \vee S \lesssim D^2 m^{14}(\log^2(1 / \eps') + \log^2(mDR\sigma^{-2})), \\
    &\log B \lesssim m^{23}( \log^3(1 / \eps') + \log^3(mDR\sigma^{-2})) .
\end{align*}
This concludes the proof.

\myendproof

\subsection{Proof of Lemma \ref{lem:Upsilon_j_gelu_approx}}
\label{sec:lem:Upsilon_j_gelu_approx_proof}

Before we move to the proof of the lemma, we state a technical result that is established in Appendix \ref{sec:lem_Psi_gelu_approx_proof}.

\begin{Lem}
\label{lem:Psi_gelu_approx}
    Assuming that Assumption \ref{asn:relax_man} holds, let $\cK_\infty$ be the compact set as defined in \eqref{eq:cK_R_def} and let $\psi : \cU_\bj \to \R$ with $\|\psi\|_{W^{0, \infty}(\cU_\bj)} \leq 2$, where $\bj \in \{1, \dots, N\}^d$ is arbitrary.
    Let also $\Psi_\bj[\psi]$ be the function given by \eqref{eq:Psi_def}.
    We also require $\eps \in (0, 1)$ be sufficiently small to satisfy
    \begin{align}
        \label{eq:choice_eps_psi_approx}
        \|g^\circ - g^*\|_{L^\infty([0, 1]^d)} \leq 1,
        \quad \frac{5 (H \vee 1)^2 P(d, \beta)^2\sqrt{D}(R + 3)\eps}{2\sigma^2} \leq 1 .
    \end{align}
    Then for any $\eps' \in (0, 1)$ and $m \in \N$ there exists a neural network $\tilde{\Psi}_\bj \in \NN(L, W, S, B)$ satisfying
    \begin{align*}
        (i) &\quad \|\tilde{\Psi}_\bj - \Psi_\bj[\psi]\|_{W^{m, \infty}(\cK_\infty)} \leq \eps^d \eps', \\
        (ii) &\quad \| \tilde{\Psi}_\bj \|_{W^{m, \infty}(\R^D)} \leq \exp\{\cO(m^3 \log(mDR \sigma^{-2}\log(1 / \eps')) )\} \eps^d .
    \end{align*}
    Moreover, $\tilde{\Psi}_\bj$ has
    \begin{align*}
        &L \asymp \log(mDR\sigma^{-2}) + \log\log(1 / \eps'),
        \quad \log B \lesssim m^{14} \log^3(1 / \eps') \log^4(mDR\sigma^{-2}), \\
        &\|W\|_\infty \vee S \lesssim D^2 ( m^3 \log(1 / \eps') + m^3\log(mDR\sigma^{-2}) )^{3 + P(d, \beta)} .
    \end{align*}
    We also note that the number of layers is the same for all $\bj \in \{1, \dots, N\}^d$.

\end{Lem}

\noindent
\textbf{Step 1: decomposing the approximation error.}\quad
Let $\tilde{\varphi}_{\bj, 0}$ denote the neural network from Lemma \ref{lem:exp_V_j_0_approx}, and $\tilde{\Psi}_\bj$ be the neural network constructed in Lemma \ref{lem:Psi_gelu_approx}, both with the precision parameter $\eps_0$.
Let also $\varphi_{mul, \bj}$ be a GELU network from Lemma \ref{lem:multi_approx_gelu} that aims to approximate the multiplication with the precision parameter $\eps_{mul} \in (0, 1)$.
The parameters $\eps_0$ and $\eps_{mul}$ will be specified in the course of the proof.
Therefore, for $\tilde{\Upsilon}_\bj = \eps^d\varphi_{mul, \bj}(\tilde{\varphi}_{\bj, 0}, \eps^{-d}\tilde{\Psi}_\bj)$, it follows that
\begin{align}
    \label{eq:tilde_ups_j_decomp}
    \notag
    &\| \tilde{\Upsilon}_\bj - \Upsilon_\bj[\psi] \|_{W^{m, \infty}(\cK_\infty)}
    \leq \underbrace{ \|\eps^d\varphi_{mul, \bj}(\tilde{\varphi}_{\bj, 0}, \eps^{-d}\tilde{\Psi}_\bj) - \tilde{\varphi}_{\bj, 0} \cdot \tilde{\Psi}_\bj \|_{W^{m, \infty}(\cK_\infty)} }_{(A)} \\
    &\quad + \underbrace{ \| (\tilde{\varphi}_{\bj, 0} - \exp \circ (-V_{\bj, 0}) ) \cdot \tilde{\Psi}_\bj \|_{W^{m, \infty}(\cK_\infty)} }_{(B)}
    + \underbrace{ \| \exp \circ (-V_{\bj, 0}) \cdot \tilde{\Psi}_\bj - \Upsilon_\bj[\psi] \|_{W^{m, \infty}(\cK_\infty)} }_{(C)} .
\end{align}
Next, we bound each term individually.

\noindent
\textbf{Stpe 2: bounding term $(A)$.}\quad
From Lemma \ref{lem:comp_sob_norm} and Lemma \ref{lem:multi_approx_gelu} we find that
\begin{align*}
    &\|\eps^d\varphi_{mul, \bj}(\tilde{\varphi}_{\bj, 0}, \eps^{-d}\tilde{\Psi}_\bj) - \eps^d\tilde{\varphi}_{\bj, 0} \cdot (\eps^{-d}\tilde{\Psi}_\bj) \|_{W^{m, \infty}(\cK_\infty)} \\
    &\quad \leq \exp\{ \cO(m\log (mD) + m\log(1 + \| \tilde{\varphi}_{\bj, 0} \|_{W^{m, \infty}(\R^D)} + \eps^{-d}\| \tilde{\Psi}_\bj \|_{W^{m, \infty}(\R^D)} )) \} \eps^d\eps_{mul}.
\end{align*}
Using claim $(ii)$ from Lemmata \ref{lem:exp_V_j_0_approx} and \ref{lem:Psi_gelu_approx}, we obtain that
\begin{align*}
    &\|\eps^d\varphi_{mul, \bj}(\tilde{\varphi}_{\bj, 0}, \eps^{-d}\tilde{\Psi}_\bj) - \tilde{\varphi}_{\bj, 0} \cdot \tilde{\Psi}_\bj \|_{W^{m, \infty}(\cK_\infty)}
    \leq \exp\{\cO( m^5\log(mDR\sigma^{-2}\log(1 / \eps_0)) )\}\eps^d\eps_{mul}.
\end{align*}
Thus, setting
\begin{align}
    \label{eq:int_eps_mul_def}
    \log(1 / \eps_{mul}) \asymp m^5\log(1 / \eps_0) + m^5\log(mDR\sigma^{-2})
\end{align}
ensures that
\begin{align}
    \label{eq:int_A_bound}
    \|\eps^d\varphi_{mul, \bj}(\tilde{\varphi}_{\bj, 0}, \eps^{-d}\tilde{\Psi}_\bj) - \eps^d\tilde{\varphi}_{\bj, 0} \cdot (\eps^{-d}\tilde{\Psi}_\bj) \|_{W^{m, \infty}(\cK_\infty)} \leq \eps^d \eps_0.
\end{align}

\noindent
\textbf{Step 3: bounding term $(B)$.}\quad
From Lemma \ref{lem:prod_sob_norm} we obtain that
\begin{align*}
    &\| (\tilde{\varphi}_{\bj, 0} - \exp \circ (-V_{\bj, 0}) ) \cdot \tilde{\Psi}_\bj \|_{W^{m, \infty}(\cK_\infty)} 
    \leq 2^m \|\tilde{\varphi}_{\bj, 0} - \exp \circ (-V_{\bj, 0}) \|_{W^{m, \infty}(\cK_\infty)} \|\tilde{\Psi}_\bj\|_{W^{m, \infty}(\cK_\infty)}.
\end{align*}
Thus, claim $(i)$ of Lemma \ref{lem:exp_V_j_0_approx} and claim $(ii)$ of Lemma \ref{lem:Psi_gelu_approx} imply that
\begin{align}
    \label{eq:int_B_bound}
    \| (\tilde{\varphi}_{\bj, 0} - \exp \circ (-V_{\bj, 0}) ) \cdot \tilde{\Psi}_\bj \|_{W^{m, \infty}(\cK_\infty)}
    \leq \exp\{\cO(m^3\log(mDR\sigma^{-2}\log(1 / \eps_0)))\} \eps^d\eps_0 .
\end{align}

\noindent
\textbf{Step 4: bounding term $(C)$.}\quad
As suggested by Lemma \ref{lem:prod_sob_norm}, it holds that
\begin{align*}
    \| \exp \circ (-V_{\bj, 0}) \cdot \tilde{\Psi}_\bj - \Upsilon_\bj[\psi] \|_{W^{m, \infty}(\cK_\infty)}
    \leq 2^m \|\exp \circ (-V_{\bj, 0})\|_{W^{m, \infty}(\cK_\infty)} \|\tilde{\Psi}_\bj - \Psi_\bj[\psi]\|_{W^{m, \infty}(\cK_\infty)}
\end{align*}
Now claim $(ii)$ of Lemma \ref{lem:exp_V_j_0_approx} and claim $(i)$ of Lemma \ref{lem:Psi_gelu_approx} yield
\begin{align}
    \label{eq:int_C_bound}
    \| \exp \circ (-V_{\bj, 0}) \cdot \tilde{\Psi}_\bj - \Upsilon_\bj[\psi] \|_{W^{m, \infty}(\cK_\infty)}
    \leq \exp\{\cO(m^4 \log(mDR\sigma^{-2} \log(1 / \eps_0)) )\} \eps^d\eps_0 .
\end{align}

\noindent
\textbf{Step 5: combining the bounds for $(A)$, $(B)$ and $(C)$.}\quad
The combination of \eqref{eq:tilde_ups_j_decomp}, \eqref{eq:int_A_bound}, \eqref{eq:int_B_bound}, and \eqref{eq:int_C_bound} immediately yields
\begin{align*}
    \| \tilde{\Upsilon}_\bj - \Upsilon_\bj[\psi] \|_{W^{m, \infty}(\cK_\infty)}
    \leq \exp\{\cO(m^4\log(mDR\sigma^{-2}\log(1 / \eps_0)))\} \eps^d\eps_0.
\end{align*}
Therefore, choosing
\begin{align}
    \label{eq:int_approx_eps_0_def}
    \log(1 / \eps_0) \asymp m^4\log(1 / \eps') + m^4\log(mDR\sigma^{-2})
\end{align}
ensures that
\begin{align*}
    \| \tilde{\Upsilon}_\bj - \Upsilon_\bj[\psi] \|_{W^{m, \infty}(\cK_\infty)}
    \leq \eps^d \eps' .
\end{align*}
Moreover, using Lemmata \ref{lem:comp_sob_norm}, \ref{lem:multi_approx_gelu}, \ref{lem:exp_V_j_0_approx}, and \ref{lem:Psi_gelu_approx}, we obtain that
\begin{align}
    \label{eq:tilde_Ups_R}
    \notag
    \|\tilde{\Upsilon}_\bj\|_{W^{m, \infty}(\R^D)}
    &\leq \exp\{\cO(m\log(mD) + m\log(\|\tilde{\varphi}_{\bj, 0}\|_{W^{m, \infty}(\R^D)} + \|\tilde{\Psi}_\bj\|_{W^{m, \infty}(\R^D)} ) )\}\eps^d \\
    &\leq \exp\{\cO(m^5 \log(mD\sigma^{-2}\log(1 / \eps')))\} \eps^d.
\end{align}

\noindent
\textbf{Step 6: deriving the configuration of $\tilde{\Upsilon}_\bj$.}\quad 
Combining \eqref{eq:int_eps_mul_def}, \eqref{eq:int_approx_eps_0_def}, and Lemma \ref{lem:multi_approx_gelu}, we conclude that $\varphi_{mul, \bj}$ belong to the neural network class $\NN(L_{mul}, W_{mul}, S_{mul}, B_{mul})$ with
\begin{align*}
    L_{mul} \vee \|W_{mul}\|_\infty \vee S_{mul} \lesssim 1,
    \quad \log B_{mul} \lesssim m^9\log(1 / \eps') + m^9\log(mDR\sigma^{-2}).
\end{align*}
From Lemmata \ref{lem:exp_V_j_0_approx} and \ref{lem:Psi_gelu_approx} we deduce that by appropriately scaling the precision parameter of either $\tilde{\varphi}_{\bj, 0}$ or $\tilde{\Psi}_\bj$,
one can ensure that the networks have the same number of layers.
Combining this observation with Lemmata \ref{lem:concat_nn}, \ref{lem:paral_nn} and the network architectures specified in
Lemmata \ref{lem:exp_V_j_0_approx} and \ref{lem:Psi_gelu_approx}, we conclude that $\tilde{\Upsilon}_\bj \in \NN(L, W, S, B)$ with
\begin{align}
    \label{eq:tilde_Upsilon_j_cfg}
    \notag
    &L \asymp \log(mDR\sigma^{-2}) + \log\log(1 / \eps') ,
    \quad \log B \lesssim m^{35}\log^3(1 / \eps')\log^7(mDR\sigma^{-2}), \\
    &\|W\|_\infty \vee S \lesssim D^2 m^{21 + 7P(d, \beta)}(\log(1 / \eps') + \log(mDR\sigma^{-2}))^{3 + P(d, \beta)} .
\end{align}
thereby finishing the proof.

\myendproof

\subsection{Proof of Lemma \ref{lem:div_num_denum_gelu_approx}}
\label{sec:lem_div_num_denum_gelu_approx_proof}

Recall from \eqref{eq:f_circ_def} that
\begin{align*}
    f^\circ_l(y) = \frac{\Upsilon_l[g^\circ]}{\Upsilon[1]},
    \quad \Upsilon_l[g^\circ] = \sum_{\bj \in \{1, \dots, N\}^d} \Upsilon_\bj[g^\circ_{\bj,l}],
    \quad \Upsilon[1] = \sum_{\bj \in \{1, \dots, N\}^d}\Upsilon_\bj[1],
    \quad 1 \leq l \leq D,
\end{align*}
where we define
\begin{align*}
    \Upsilon_\bj[g^\circ_{\bj, l}] = \integral{\cU_\bj}g^\circ_{\bj, l}(u) \exp\left\{-\frac{\|y - g^\circ_\bj(u)\|^2}{ 2\sigma^2} \right\} \dd u,
    \quad 
    \Upsilon_\bj[1] = \integral{\cU_\bj} \exp\left\{-\frac{\|y - g^\circ_\bj(u)\|^2}{ 2\sigma^2} \right\} \dd u .
\end{align*}

\noindent
\textbf{Step 1: decomposing the approximation error.}\quad
Lemma \ref{lem:div_gelu_approx} yields that for each $1 \leq l \leq D$ there exists a division GELU network $\varphi_{div, l}$ with the precision parameter $\eps'$ and the scale parameter $N_{div} \in \N$ with $N_{div} \geq 3$, both of which will be optimized later in the proof.
We also put the smoothness parameter $m + 1$.
In addition, by rescaling the inputs of each $\varphi_{div, l}$ by a constant factor, we generalize the approximation guarantee to
\begin{align}
    \label{eq:phi_div_l_acc}
    \max_{1 \leq l \leq D}\|\varphi_{div, l} - \mathrm{div}\|_{W^{m, \infty}([-3, 3] \times [2^{-N - 1}, 2])} \leq \eps'.
\end{align}
Thus, for all $1 \leq l \leq D$ the triangle inequality suggests that
\begin{align}
    \label{eq:phi_div_l_A_B_decomp}
    \notag
    &\|\varphi_{div, l}(\cP_l, \cQ) - f^\circ_l\|_{W^{m, \infty}(\cK_\infty)} \\
    &\quad \leq \underbrace{ \|\varphi_{div, l}(\cP_l, \cQ) - \varphi_{div, l}(\Upsilon_l[g^\circ], \Upsilon[1])\|_{W^{m, \infty}(\cK_\infty)} }_{(A)}
    + \underbrace{ \| \varphi_{div, l}(\Upsilon_l[g^\circ], \Upsilon[1]) - f^\circ_l \|_{W^{m, \infty}(\cK_\infty)} }_{(B)} .
\end{align}
We now establish a lower bound for the denominator $\Upsilon[1]$.
To this end, recall from the definition of $\cK$ given in \eqref{eq:cK_R_def} that 
\begin{align}
\label{eq:div_R_inf_bound}
    R_\infty = \frac{5}{2}\vee \sup_{y \in \cK} \|y\|_\infty
    \leq \frac{5}{2} \vee \sup_{y \in \cK} \inf_{u \in [0, 1]^d}(\|y - g^*(u)\| + \|g^*(u)\|)
    \leq R + 4 .
\end{align}
Consequently, for any $y \in \cK_\infty$, we have that
\begin{align*}
    \sum_{\bj \in \{1, \dots, N\}^d }\integral{\cU_\bj}\exp\left\{-\frac{\|y - g^\circ_\bj(u)\|^2}{2\sigma^2}\right\} \dd u
    &\geq \exp\left\{ -\frac{(R + 4)^2}{\sigma^2} - \frac{4}{\sigma^2} \right\}
    \geq \exp\left\{ -\frac{2R^2}{\sigma^2} - \frac{36}{\sigma^2} \right\} .
\end{align*}
In view of the derived lower bound, we choose the scale parameter for $\varphi_{div, l}$ for each $1 \leq l \leq D$ such that
\begin{align}
    \label{eq:N_def}
    N_{div} \asymp (R^2 + 1) / \sigma^2 ,
\end{align}
thereby guaranteeing that
\begin{align}
    \label{eq:Upsilon_1_lb_N}
    \min_{y \in \cK}\Upsilon[1](y) \geq 2^{-N_{div}}.
\end{align}

\noindent
\textbf{Step 2: bounding term $(A)$.}\quad
Assume that the parameter $\eps' \leq 2^{-N_{div} - 1}$.
Hence, from \eqref{eq:p_q_networks_acc_K}, \eqref{eq:Upsilon_1_lb_N} and Lemma \ref{lem:comp_sob_norm} it follows that
\begin{align*}
    &\|\varphi_{div, l}(\cP_l, \cQ) - \varphi_{div, l}(\Upsilon_l[g^\circ], \Upsilon[1])\|_{W^{m, \infty}(\cK_\infty)}
    \leq \|\varphi_{div, l}\|_{W^{m + 1, \infty}([-3, 3] \times [2^{-N_{div}-1}, 2])}\eps' \\
    &\quad \cdot \exp\{ m\log(\|\cP_l\|_{W^{m, \infty}(\R^D)} + \|\cQ\|_{W^{m, \infty}(\R^D)}) )\}  .
\end{align*}
Therefore, using \eqref{eq:Q_P_l_bound_R}, \eqref{eq:phi_div_l_acc} and Lemma \ref{lem:div_gelu_approx}, we obtain that
\begin{align*}
    &\|\varphi_{div, l}(\cP_l, \cQ) - \varphi_{div, l}(\Upsilon_l[g^\circ], \Upsilon[1])\|_{W^{m, \infty}(\cK_\infty)} \\
    &\quad \leq \exp\{\cO(m^6 \log(mDR\sigma^{-2}\log(1 / \eps')))\} \|\mathrm{div}\|_{W^{m, \infty}([-3, 3] \times [2^{-N_{div}-1}, 2])} \eps' .
\end{align*}
Noting that
\begin{align}
    \label{eq:div_three_2_N_bound}
    \|\mathrm{div}\|_{W^{m, \infty}([-3, 3] \times [2^{-N_{div}-1}, 2])} \leq \exp\{\cO(m\log m + mN_{div})\},
\end{align}
we conclude that
\begin{align}
    \label{eq:phi_div_l_A_bound}
    &\|\varphi_{div, l}(\cP_l, \cQ) - \varphi_{div, l}(\Upsilon_l[g^\circ], \Upsilon[1])\|_{W^{m, \infty}(\cK_\infty)}
    \leq \exp\{\cO(m^6 N_{div} \log(mDR\sigma^{-2}\log(1 / \eps')))\} \eps' .
\end{align}

\noindent
\textbf{Step 3: bounding term $(B)$.}\quad
We find that \eqref{eq:phi_div_l_acc} together with Lemma \ref{lem:comp_sob_norm} implies that
\begin{align*}
    &\| \varphi_{div, l}(\Upsilon_l[g^\circ], \Upsilon[1]) - f^\circ_l \|_{W^{m, \infty}(\cK_\infty)} \\
    &\quad \leq \exp\{\cO(m\log(mD + \|\Upsilon_l[g^\circ]\|_{W^{m, \infty}(\cK_\infty)} + \|\Upsilon[1] \|_{W^{m, \infty}(\cK_\infty)} ))\} \eps' .
\end{align*}
Applying \eqref{eq:p_q_networks_acc_K}, \eqref{eq:Q_P_l_bound_R} and \eqref{eq:phi_div_l_acc}, we obtain that
\begin{align}
    \label{eq:phi_div_l_B_bound}
    \| \varphi_{div, l}(\Upsilon_l[g^\circ], \Upsilon[1]) - f^\circ_l \|_{W^{m, \infty}(\cK_\infty)}
    \leq \exp\{\cO(m^6 \log(mDR\sigma^{-2}\log(1 / \eps')))\} \eps' .
\end{align}

\noindent
\textbf{Step 4: combining the bounds for terms $(A)$ and $(B)$.}\quad
Substituting the bounds \eqref{eq:phi_div_l_A_bound} and \eqref{eq:phi_div_l_B_bound} into \eqref{eq:phi_div_l_A_B_decomp}, we obtain that
\begin{align*}
    \|\varphi_{div, l}(\cP_l, \cQ) - f^\circ_l\|_{W^{m, \infty}(\cK_\infty)}
    \leq \exp\{\cO(m^6 N_{div} \log(mDR\sigma^{-2}\log(1 / \eps')))\} \eps' .
\end{align*}
Therefore, choosing
\begin{align}
    \label{eq:eps_prime_def}
    \log(1 / \eps')
    \asymp m^6 N_{div} \log(mDR\sigma^{-2})\log(1 / \eps)
    \asymp m^6 (R^2 + 1)\sigma^{-2} \log(mDR\sigma^{-2})\log(1 / \eps) ,
\end{align}
as suggested by \eqref{eq:N_def}, ensures that $\eps' \leq 2^{-N_{div} - 1}$ and
\begin{align*}
    \|\varphi_{div, l}(\cP_l, \cQ) - f^\circ_l\|_{W^{m, \infty}(\cK_\infty)}
    \leq \eps^\beta / 2.
\end{align*}
According to Lemma \ref{lem:paral_nn}, a parallel stack of $\{\varphi_{div, l}(\cP_l, \cQ)\}_{l=1}^D$ denoted as $\tilde{f}$ satisfies
\begin{align}
    \label{eq:tilde_f_l_f_circ_acc}
    \max_{1 \leq l \leq D}\|\tilde{f}_l - f^\circ_l\|_{W^{m, \infty}(\cK_\infty)} \leq \eps^\beta / 2.
\end{align}
Moreover, from \eqref{eq:Q_P_l_bound_R}, Lemma \ref{lem:comp_sob_norm} and Lemma \ref{lem:div_gelu_approx} we find that
\begin{align*}
    \max_{1 \leq l \leq D}\|\tilde{f}_l\|_{W^{m, \infty}(\R^D)}
    &\leq \exp\{\cO(m\log(mD + \|\cP_l\|_{W^{m, \infty}(\R^D)} + \|\cQ\|_{W^{m, \infty}(\R^D)}))\} \|\varphi_{div, l}\|_{W^{m, \infty}(\R^2)} \\
    &\leq \exp\{\cO(m^6 \log(mDR\sigma^{-2}\log(1 / \eps')) + m^4 N_{div} + m^4\log(m\log(1 / \eps')))\} .
\end{align*}
Hence, substituting the values of $N_{div}$ from \eqref{eq:N_def} and $\eps'$ from \eqref{eq:eps_prime_def} yields
\begin{align*}
    \max_{1 \leq l \leq D}\|\tilde{f}_l\|_{W^{m, \infty}(\R^D)}
    \leq \exp\{ \cO( m^6(R^2 + 1)\sigma^{-2}\log(mDR\sigma^{-2}\log(1 / \eps)) )\} .
\end{align*}
From Lemma \ref{lem:f_circ_f_star_acc} we deduce that
\begin{align*}
    \max_{1 \leq l \leq D} (\|f^\circ_l\|_{W^{m, \infty}(\R^D)} \vee \|f^*_l\|_{W^{m, \infty}(\R^D)} ) \leq \exp\{\cO(m\log(m\sigma^{-2}))\} .
\end{align*}
Hence, this observation together with \eqref{eq:tilde_f_l_f_circ_acc} implies that
\begin{align}
    \label{eq:tilde_f_bound_R}
    \max_{1 \leq l \leq D} \|\tilde{f}_l\|_{W^{m, \infty}(\cK_\infty)}
    \leq \exp\{\cO(m\log(m\sigma^{-2}))\} .
\end{align}

\noindent
\textbf{Step 5: clipping the input.}\quad
Now let $\{ \tilde{\psi}_{clip, l} \}_{l = 1}^D$ be the clipping function approximators from Lemma \ref{lem:clip_gelu_approx} with the precision parameter $\eps_{\psi, clip}$ and the scale parameter $R_\infty$ as defined in \eqref{eq:K_infty_def}.
Applying Lemma \ref{lem:paral_nn}, we conclude that there exists a neural network $\psi_{clip}$ that performs component-wise clipping of a $D$-dimensional vector.
Formally,
\begin{align*}
    \max_{1 \leq l \leq D} \|(\psi_{clip})_l - (\id)_l\|_{W^{m, \infty}(\cK_\infty)} \leq \eps_{\psi, clip} .
\end{align*}
Therefore, Lemma \ref{lem:comp_sob_norm} suggests that for $\cK_{\infty, clip} = \cK_\infty \cup \psi_{clip}(\cK_\infty)$, we have that
\begin{align*}
    \max_{1 \leq l \leq D}\| \tilde{f}_l \circ \psi_{clip} - \tilde{f}_l \|_{W^{m, \infty}(\cK_\infty)}
    \leq \exp\{\cO(m\log(m D + R_\infty) )\} \| \tilde{f}_l \|_{W^{m + 1, \infty}(\cK_{\infty, clip} )} \eps_{\psi, clip} .
\end{align*}
Using the bound given in \eqref{eq:div_R_inf_bound} and the observation that $\psi_{clip}(\cK_{\infty}) \in 2\cK_\infty$ due to Lemma \ref{lem:clip_gelu_approx}, we find that 
\begin{align*}
    \max_{1 \leq l \leq D}\| \tilde{f}_l \circ \psi_{clip} - \tilde{f}_l \|_{W^{m, \infty}(\cK_\infty)}
    \leq \exp\{ m\log(mDR) \} \| \tilde{f}_l \|_{W^{m + 1, \infty}(2 \cK_\infty )} \eps_{\psi, clip} .
\end{align*}
We emphasize that the approximation result \eqref{eq:tilde_f_l_f_circ_acc} can be extended to the compact $2\cK_\infty$ and the smoothness parameter $m + 1$,
by only modifying the configuration of the resulting neural network by an absolute constant.
Thus, using \eqref{eq:tilde_f_bound_R}, we deduce that
\begin{align*}
    \max_{1 \leq l \leq D}\| \tilde{f}_l \circ \psi_{clip} - \tilde{f}_l \|_{W^{m, \infty}(\cK_\infty)}
    \leq \exp\{ \cO(m\log(mDR\sigma^{-2})) \} \eps_{\psi, clip} .
\end{align*}
Hence, choosing
\begin{align}
    \label{eq:eps_psi_clip_def}
    \log(1 / \eps_{\psi, clip}) \asymp m\log(mDR\sigma^{-2}) + \log(1 / \eps)
\end{align}
ensures that
\begin{align*}
    \max_{1 \leq l \leq D}\| \tilde{f}_l \circ \psi_{clip} - \tilde{f}_l \|_{W^{m, \infty}(\cK_\infty)}
    \leq \eps^\beta / 2.
\end{align*}
The combination of the derived bound and \eqref{eq:tilde_f_l_f_circ_acc} yields
\begin{align*}
    \max_{1 \leq l \leq D}\|\tilde{f}_l \circ \psi_{clip} - f^\circ_l\|_{W^{m, \infty}(\cK_\infty)} \leq \eps^\beta.
\end{align*}
Furthermore, from Lemmata \ref{lem:comp_sob_norm} and \ref{lem:clip_gelu_approx} we deduce that for any $1 \leq k \leq m$ it holds that
\begin{align*}
    |\tilde{f}_l \circ \psi_{clip}|_{W^{k, \infty}(\R^D)}
    &\leq \exp\{\cO(k \log(kD))\}|\tilde{f}_l|_{W^{k, \infty}(2\cK_\infty)} (1 \vee \|\psi_{clip}\|_{W^{k, \infty}(\R^D)})^k \\
    &\leq \exp\{\cO(k \log(kD) + k^2\log m + k^2\log\log(1 / \eps_{\psi, clip}))\}(1 + |f^\circ_l|_{W^{k, \infty}(2\cK_\infty)}) ,
\end{align*}
where the last inequality uses \eqref{eq:tilde_f_l_f_circ_acc}.
Hence, Lemma \ref{lem:f_circ_f_star_acc} and the choice of $\eps_{\psi, clip}$ given in \eqref{eq:eps_psi_clip_def} imply that
\begin{align*}
    |\tilde{f}_l \circ \psi_{clip}|_{W^{k, \infty}(\R^D)}
    \leq \sigma^{-2k} \exp\{\cO( k^2\log(mD\log(1 / \eps)\log(R\sigma^{-2})) )\}, \quad \text{for all } 1 \leq k \leq m .
\end{align*}
As for the case $k = 0$, from \eqref{eq:tilde_f_l_f_circ_acc} and Lemma \ref{lem:f_circ_f_star_acc} it follows that
\begin{align*}
    |\tilde{f}_l \circ \psi_{clip}|_{W^{0, \infty}(\R^D)}
    \leq |\tilde{f}_l|_{W^{0, \infty}(2\cK_\infty)}
    \leq 1 + |f^\circ_l|_{W^{0, \infty}(2\cK_\infty)}
    \leq \cO(1) .
\end{align*}

\noindent
\textbf{Step 6: deriving the configuration.}\quad
Next, we derive the configuration of $\bar{f} = \tilde{f} \circ \psi_{clip}$.
We note from \eqref{eq:eps_psi_clip_def} and Lemma \ref{lem:clip_gelu_approx} that $\psi_{clip} \in \NN(L_\psi, W_\psi, S_\psi, B_\psi)$ with
\begin{align*}
    L_\psi \lesssim 1,
    \quad \|W_\psi\|_\infty \vee S_\psi \lesssim D,
    \quad \log B_\psi \lesssim m\log(mDR\sigma^{-2}) + \log(1 / \eps) .
\end{align*}
In addition, from \eqref{eq:N_def}, \eqref{eq:eps_prime_def} and Lemma \ref{lem:div_gelu_approx} we find that for each $1 \leq l \leq D$ the GELU network $\varphi_{div, l} \in \NN(L_{div}, W_{div}, S_{div}, B_{div})$ with
\begin{align*}
    &L_{div} \lesssim \log(mDR\sigma^{-2}\log(1 / \eps)),
    \quad \|W_{div}\|_\infty \vee S_{div} \lesssim m^{45}R^{18}\sigma^{-18} \log^4(mDR\sigma^{-2}) \log^4(1 / \eps), \\
    &\log B_{div} \lesssim  m^{48}R^{16}\sigma^{-16} \log^4(mDR\sigma^{-2}) \log^{4}(1 / \eps) .
\end{align*}
Therefore, combining the derived configurations with those for the networks $\cQ, \cP_1, \dots, \cP_D$ given in \eqref{eq:cQ_cP_cfg},
and using the choice of $\eps'$ from \eqref{eq:eps_prime_def} together with Lemmata \ref{lem:concat_nn} and \ref{lem:paral_nn}, we conclude that $\bar{f} \in \NN(\bar{L}, \bar{W}, \bar{S}, \bar{B})$ with
\begin{align*}
    &\bar{L} \lesssim \log(mDR\sigma^{-2}\log(1 / \eps)),
    \quad \log \bar{B} \lesssim m^{53}R^{16}\sigma^{-16}\log^{10}(mDR\sigma^{-2})\log^4(1 / \eps), \\
    &\|\bar{W}\|_\infty \vee \bar{S}
    \lesssim \eps^{-d} D^4 m^{84 + 13 P(d, \beta)} (R\sigma^{-2})^{24 + 2 P(d, \beta)} (\log(1 / \eps) \log(mDR\sigma^{-2}))^{7 + P(d,\beta)} .
\end{align*}
Thus, the proof is complete.

\myendproof

\subsection{Proof of Lemma \ref{lem:p_star_tail}}
\label{sec:lem_p_star_tail_proof}
The proof follows the identical argument as provided in Lemma A.2 of \cite{yakovlev2025generalization}.
The definition of $\sfp^*$ given in \eqref{eq:p_star_def} yields
\begin{align*}
    \integral{\R^D \setminus \cK}\sfp^*(y) \dd y
    = \integral{[0, 1]^d}\integral{\R^D \setminus \cK} (\sqrt{2\pi}\sigma)^{-D} \exp\left\{-\frac{\|y - g^*(u)\|}{2\sigma^2}\right\} \dd y \, \dd u .
\end{align*}
Now for an arbitrary $u \in [0, 1]^d$ we introduce a random variable $Y \sim \cN(g^*(u), \sigma^2 I_D)$.
Therefore, the inner integral is evaluated as
\begin{align*}
    \integral{\R^D \setminus \cK} (\sqrt{2\pi}\sigma)^{-D} \exp\left\{-\frac{\|y - g^*(u)\|}{2\sigma^2}\right\} \dd y
    \leq \p(\|Y - g^*(u)\| \geq R)
    = \p\left(\frac{\|Y - g^*(u)\|}{\sigma^2} \geq \frac{R^2}{\sigma^2}\right) ,
\end{align*}
where $\sigma^{-2}\|Y - g^*(u)\|^2 \sim \chi^2(D)$.
Invoking the concentraion bound for chi-squared distribution outlined in Lemma \ref{lem:chi_sq_concentr}, we deduce that
\begin{align*}
    \p\left(\frac{\|Y - g^*(u)\|}{\sigma^2} \geq \frac{R^2}{\sigma^2}\right)
    &= \p\left(\frac{\|Y - g^*(u)\|}{\sigma^2} \geq D + \frac{R^2 - \sigma^2 D}{\sigma^2}\right) \\
    &\leq \exp\left\{-\frac{1}{16}\left(\frac{R^2 - \sigma^2 D}{\sigma^2 D} \wedge \frac{\sqrt{R^2 - \sigma^2 D}}{\sigma}\right)\right\} .
\end{align*}
Therefore, we conclude that
\begin{align*}
    \integral{\R^D \setminus \cK} \sfp^*(y) \dd y
    \leq \exp\left\{-\frac{1}{16}\left(\frac{R^2 - \sigma^2 D}{\sigma^2 D} \wedge \frac{\sqrt{R^2 - \sigma^2 D}}{\sigma}\right)\right\} ,
\end{align*}
thereby completing the proof.

\myendproof

\subsection{Proof of Lemma \ref{lem:Psi_gelu_approx}}
\label{sec:lem_Psi_gelu_approx_proof}
The proof is quite technical.
To this end, we split it into several steps.

\noindent
\textbf{Step 1: building an approximation on a compact set.}\quad
First, use Taylor's expansion of the exponent and analyze its accuracy in the Sobolev norm.
Formally, for some $r \in \N$, which will be adjusted later in the proof, let
\begin{align*}
    \overline{\exp}(x) = \sum_{i=0}^{r - 1}\frac{x^i}{i!}, \quad x \in \R.
\end{align*}
Therefore, for all $0 \leq k \leq m$ and $r > k$ Stirling's approximation yields
\begin{align*}
    |\exp - \overline{\exp}|_{W^{k, \infty}([-1, 1])} \leq \frac{e}{(r - k)!}
    \leq e \left(\frac{e}{r - k}\right)^{r - k}.
\end{align*}
Hence, setting $r = \ceil{m + 1 + e^2 + \log(1 / \eps_0)}$ for some $\eps_0 \in (0, 1)$, which will be optimized later, we obtain that
\begin{align}
    \label{eq:exp_bar_acc_ints}
    \|\exp - \overline{\exp}\|_{W^{m, \infty}([-1, 1])} \leq \eps_0 .
\end{align}
The triangle inequality together with the fact that $\|\psi\|_{W^{0, \infty}(\cU_\bj)} \leq 2$ suggests that
\begin{align}
    \label{eq:Psi_j_bar_exp_acc}
    \notag
    &\left\|\Psi_\bj[\psi] - \integral{\cU_\bj}\psi(u) \overline{\exp}(-\sR_\bj^\top a_\bj(u)) \dd u \right\|_{W^{m, \infty}(\cK_\infty)} \\
    \notag
    &\quad \leq \integral{\cU_\bj}|\psi(u)| \cdot \left\| (\exp - \overline{\exp}) \circ (-\sR_\bj^\top a_\bj(u)) \right\|_{W^{m, \infty}(\cK_\infty)} \dd u \\
    &\quad \leq 2 \, \mathrm{vol}(\cU_\bj) \sup_{u \in \cU_\bj} \left\| (\exp - \overline{\exp}) \circ (-\sR_\bj^\top a_\bj(u)) \right\|_{W^{m, \infty}(\cK_\infty)}.
\end{align}
Next, for an arbitrary $u \in \cU_\bj$ we evaluate
\begin{align*}
    &\|\sR_\bj^\top a_\bj(u)\|_{W^{m, \infty}(\cK_\infty)} \\
    &\quad \leq \frac{\|g^\circ_\bj(u) - g^\circ_\bj(u_\bj)\|^2}{2\sigma^2} + P(d, \beta) \max_{\bk \in \Z_+^d, \; 1 \leq |\bk| \leq \floor{\beta}} \frac{(u - u_\bj)^\bk}{\sigma^2 \bk !} \|(\cdot - g^\circ_\bj(u_\bj))^\top \partial^\bk g^*(u_\bj)\|_{W^{m, \infty}(\cK_\infty)}.
\end{align*}
Since Assumption \ref{asn:relax_man} is fulfilled, we deduce that
\begin{align}
    \label{eq:g_circ_u_g_circ_uj_acc}
    \sup_{u \in \cU_\bj}\|g^\circ_\bj(u) - g^\circ_\bj(u_\bj)\|
    = \sup_{u \in \cU_\bj}\left\|\sum_{\bk \in \Z_+^d, \; 1 \leq |\bk| \leq \floor{\beta}}\frac{\partial^\bk g^*(u_\bj) (u - u_\bj)^\bk}{\bk!}\right\|
    \leq H \cdot P(d, \beta) \eps .
\end{align}
We also note from \eqref{eq:K_infty_def} that
\begin{align}
    \label{eq:sup_y_k_inf}
    \sup_{y \in \cK_\infty}\|y\| 
    \leq (5/2 + \sup_{y \in \cK}\|y\|)\sqrt{D}
    \leq (R + 3)\sqrt{D} .
\end{align}
Therefore, we conclude that
\begin{align}
    \label{eq:sup_u_r_a_comp}
    \notag
    \sup_{u \in \cU_\bj}\|\sR_\bj^\top a_\bj(u)\|_{W^{m, \infty}(\cK_\infty)}
    &\leq \frac{P(d, \beta)^2 H^2 \eps^2}{2\sigma^2} + \frac{P(d, \beta) H (\|g^\circ\|_{L^\infty([0, 1]^d)} + \sup_{y \in \cK_\infty}\|y\|) \eps}{\sigma^2} \\
    &\leq \frac{5 (H \vee 1)^2 P(d, \beta)^2\sqrt{D}(R + 3)\eps}{2\sigma^2} .
\end{align}
The choice of $\eps$ guarantees that the resulting bound is at most one.
Thus, applying Lemma \ref{lem:comp_sob_norm}, we obtain that for any $u \in \cU_\bj$
\begin{align*}
    &\|(\exp - \overline{\exp}) \circ (-\sR_\bj^\top a_\bj(u)) \|_{W^{m, \infty}(\cK_\infty)} \\
    &\quad \leq \exp\{\cO(m\log(mD + \|\sR_\bj^\top a_\bj(u)\|_{W^{m, \infty}(\cK_\infty)} ))\} \|\exp - \overline{\exp}\|_{W^{m, \infty}([-1, 1])} .
\end{align*}
Consequently, the combination of \eqref{eq:exp_bar_acc_ints} and \eqref{eq:sup_u_r_a_comp} leads to
\begin{align*}
    &\|(\exp - \overline{\exp}) \circ (-\sR_\bj^\top a_\bj(u)) \|_{W^{m, \infty}(\cK_\infty)}
    \leq \exp\{\cO(m\log(mD))\}\eps_0 .
\end{align*}
Substituting the derived bound into \eqref{eq:Psi_j_bar_exp_acc} implies that
\begin{align}
    \label{eq:Psi_int_exp_over_acc}
    \left\|\Psi_\bj[\psi] - \integral{\cU_\bj}\psi(u) \overline{\exp}(-\sR_\bj^\top a_\bj(u)) \dd u \right\|_{W^{m, \infty}(\cK_\infty)}
    \leq \exp\{\cO(m \log(m D))\} \eps^d \eps_0 .
\end{align}
We now integrate out the variable $u \in \cU_\bj$ and obtain that for any $y \in \R^D$ the following holds
\begin{align*}
    \integral{\cU_\bj}\psi(u)\overline{\exp}(-\sR_\bj(y)^\top a_\bj(u)) \dd u
    &= \sum_{i = 0}^{r - 1}\frac{(-1)^i}{i!}\integral{\cU_\bj} \psi(u) (\sR_\bj(y)^\top a_\bj(u) )^i \dd u \\
    &= \sum_{i = 0}^{r - 1}\frac{(-1)^i}{i!}\sum_{\substack{\bk \in \Z_+^{P(d, \beta)}, \; |\bk| = i}} \frac{i!}{\bk!} \sR_\bj(y)^\bk \integral{\cU_\bj}\psi(u)a_\bj^\bk(u) \dd u,
\end{align*}
where the last equality uses the multinomial theorem.
Introduce
\begin{align*}
    a_{\bj, i, \bk} = \frac{(-1)^i}{\bk!} \integral{\cU_\bj}\psi(u)a_\bj^\bk(u) \dd u,
    \quad 0\leq i \leq r - 1, \; \bk \in \Z_+^{P(d, \beta)}, \; |\bk| = i.
\end{align*}
Then, the approximation takes the following form:
\begin{align*}
    \integral{\cU_\bj}\psi(u)\overline{\exp}(-\sR_\bj(y)^\top a_\bj(u)) \dd u
    = \sum_{i = 0}^{r - 1} \sum_{\bk \in \Z_+^{P(d, \beta)}, \; |\bk| = i} a_{\bj, i, \bk} \cdot \sR_\bj(y)^\bk.
\end{align*}
In addition, we note from \eqref{eq:a_j_def}, \eqref{eq:g_circ_u_g_circ_uj_acc} and the choice of $\eps$ given by \eqref{eq:choice_eps_psi_approx} that
\begin{align}
    \label{eq:a_j_infty_bound}
    \sup_{u \in \cU_\bj}\|a_\bj(u)\|_\infty \leq \eps \vee H^2 P(d, \beta)^2 \eps^2 \leq \eps,
\end{align}
which leads to
\begin{align*}
    \max_{0\leq i \leq r - 1, \; \bk \in \Z_+^{P(d, \beta)}, \; |\bk| = i}|a_{\bj, i, \bk}|
    \leq \mathrm{Vol}(\cU_\bj) \cdot \|\psi\|_{W^{0, \infty}(\cU_\bj)} (\sup_{u \in \cU_\bj}\|a_\bj(u)\|_{\infty} \vee 1)^r
    \leq 2\eps^d,
\end{align*}
It is clear from \eqref{eq:r_j_def} that $\sR_\bj$ admits exact implementation with a single-layer GELU network.
We also find from \eqref{eq:v_j_0_k_def} and \eqref{eq:sup_y_k_inf} that
\begin{align}
    \label{eq:r_j_sob_bound}
    \max_{1 \leq i \leq P(d, \beta)}\|(\sR_\bj)_i\|_{W^{m, \infty}(\cK_\infty)}
    \lesssim \frac{1 + \sup_{y \in \cK_\infty} \|y\|}{\sigma^2}
    \lesssim \frac{R \sqrt{D}}{\sigma^2} .
\end{align}
Therefore, applying Lemma \ref{lem:part_sum_monomial} with $I = P(d, \beta)$ and $d = r - 1$, the accuracy parameter $\eps_{\bj, sum} \in (0, \eps]$ and the scale parameter given in \eqref{eq:r_j_sob_bound},
we obtain that there exists $\varphi_{\bj, sum} \in \NN(L_{sum}, W_{sum}, S_{sum}, B_{sum})$ such that
\begin{align}
    \label{eq:phi_j_sum_bound}
    \left\|\varphi_{\bj, sum} - \sum_{i = 0}^{r - 1}\frac{a_{\bj, i, \bk}}{2\eps^d} \cdot \mathrm{prod}_{\bk} \right\|_{W^{m, \infty}(\sR_\bj(\cK_\infty))}
    \leq \eps_{\bj, sum}.
\end{align}
Furthermore, $\varphi_{\bj, sum}$ has
\begin{align}
    \label{eq:phi_j_sum_cfg}
    \notag
    &L_{sum} \lesssim \log r,
    \quad \|W_{sum}\|_\infty \vee S_{sum} \lesssim r^{3 + P(d, \beta)}, \\
    &\log B_{sum} \lesssim (\log(1 / \eps_{\bj, sum}) + m^2 r\log(mr D R \sigma^{-2}) + m^2r^2 ) \log r .
\end{align}
Hence, the approximation accuracy for $\breve{\Psi}_\bj = 2\eps^d \cdot \varphi_{\bj, sum} \circ \sR_\bj$ follows from Lemma \ref{lem:comp_sob_norm}:
\begin{align*}
    &\left\|\breve{\Psi}_\bj - \integral{\cU_\bj}\psi(u)\overline{\exp}(-\sR_\bj^\top a_\bj(u)) \dd u \right\|_{W^{m, \infty}(\cK_\infty)} \\
    &\leq \exp\{\cO(m\log(m D + \max_{1 \leq i \leq \dim a_\bj}\|(\sR_\bj)_i\|_{W^{m, \infty}(\cK_\infty)} ) )\}\eps^d \left\|\varphi_{\bj, sum} - \sum_{i = 0}^{r - 1}\frac{a_{\bj, i, \bk}}{2\eps^d} \cdot \mathrm{prod}_{\bk} \right\|_{W^{m, \infty}(\sR_\bj(\cK_\infty))} ,
\end{align*}
Consequently, due to \eqref{eq:r_j_sob_bound} and \eqref{eq:phi_j_sum_bound}, the approximation error is bounded by
\begin{align*}
    \left\|\breve{\Psi}_\bj - \integral{\cU_\bj}\psi(u)\overline{\exp}(-\sR_\bj^\top a_\bj(u)) \dd u \right\|_{W^{m, \infty}(\cK_\infty)}
    &\leq \exp\{\cO( m\log(m D R \sigma^{-2}) )\}\eps^d \eps_{\bj, sum} .
\end{align*}
Therefore, setting
\begin{align}
    \label{eq:eps_b_sum_def}
    \log(1 / \eps_{\bj, sum})
    \asymp \log(1 / \eps_0) + m\log(mD R \sigma^{-2})
\end{align}
guarantees that
\begin{align}
    \label{eq:tilde_ups_aprox}
    \left\|\breve{\Psi}_\bj - \integral{\cU_\bj}\psi(u)\overline{\exp}(-\sR_\bj^\top a_\bj(u)) \dd u \right\|_{W^{m, \infty}(\cK_\infty)}
    \leq \eps^d \eps_0.
\end{align}
In addition, the configuration of $\breve{\Psi}_\bj \in \NN(\breve{L}, \breve{W}, \breve{S}, \breve{B})$, in view of \eqref{eq:phi_j_sum_cfg}, \eqref{eq:eps_b_sum_def} and Lemma \ref{lem:concat_nn}, takes the following form:
\begin{align}
    \label{eq:cfg_tilde_ups}
    \notag
    &\breve{L} \lesssim \log(m + \log(1 / \eps_0)),
    \quad \|\breve{W}\|_\infty \vee \breve{S} \lesssim D (m + \log(1 / \eps_0))^{3 + P(d, \beta)} \\
    &\log \breve{B} \lesssim m^5 \log^3(1 / \eps_0) \log(mDR\sigma^{-2}) ,
\end{align}
where we used the fact that $r \lesssim m + \log(1 / \eps_0)$.

\noindent
\textbf{Step 2: clipping the network input.}\quad
To make the approximation properly bounded at the infinity, we use a clipping operation.
Formally, let $\varphi_{clip}$ be a clipping operation from Lemma \ref{lem:clip_gelu_approx} with the accuracy parameter $\eps_0$ and the clipping parameter given by \eqref{eq:sup_y_k_inf}.
Let $\varphi_{\bj, clip}$ be a parallelization of $D$ independent copies of $\varphi_{clip}$ that approximates a component-wise clipping of a $D$-dimensional vector.
Then, it holds that the approximation accuracy is 
\begin{align}
    \label{eq:phi_j_clip_acc}
    \max_{1 \leq i \leq D}\|(\varphi_{\bj, clip})_i - (\id)_i\|_{W^{m, \infty}(\cK_\infty)}
    \leq \eps_0
\end{align}
and also
\begin{align}
    \label{eq:phi_j_clip_R}
    \notag
    \max_{1 \leq i \leq D }\|(\varphi_{\bj, clip})_i\|_{W^{m, \infty}(\R^D)}
    &\leq \exp\{\cO(m\log (mD) + m\log\log(1 / \eps_0) + \log(DR\sigma^{-2}))\} \\
    &\leq \exp\{\cO(m \log(mDR \sigma^{-2}) + m\log\log(1 / \eps_0) )\} .
\end{align}
Statement $(iv)$ from Lemma \ref{lem:clip_gelu_approx} in conjunction with the definition of $\cK_\infty$ outlined in \eqref{eq:K_infty_def} guarantees that $\cK_\bj := \cK_\infty \cup \varphi_{\bj, clip}(\R^{P(d, \beta)}) \subseteq 2\cK_\infty$.
Hence, one can prove a similar bound as in \eqref{eq:sup_u_r_a_comp} for $2\cK_\infty$.
Consequently, this observation together with Lemma \ref{lem:comp_sob_norm} yields 
\begin{align*}
    &\|\Psi_\bj[\psi]\|_{W^{m, \infty}(\cK_\bj)} \\
    &\quad \leq \exp\{\cO(m \log(mD + \max_{1 \leq i \leq \dim P(d, \beta)}\|(\sR_\bj)_i\|_{W^{m, \infty}(\cK_\bj)} ))\}\left\|\integral{\cU_\bj} \psi(u) \exp(-\sR_\bj^\top a_\bj(u)) \dd u \right\|_{W^{m, \infty}(\sR_\bj(\cK_\bj))} \\
    &\quad \leq \exp\{\cO(m \log(m D R \sigma^{-2}))\}\eps^d \exp\{\sup_{u \in \cU_\bj}\|\sR_\bj^\top a_\bj(u)\|_{W^{m, \infty}(2\cK_\infty)} \} .
\end{align*}
The choice of $\eps$ given in \eqref{eq:choice_eps_psi_approx} implies that
\begin{align}
    \label{eq:Psi_kj_bound}
    \|\Psi_\bj[\psi]\|_{W^{m, \infty}(\cK_\bj)} \leq \exp\{\cO(m \log(m D R \sigma^{-2}))\}\eps^d .
\end{align}
Note that the bounds \eqref{eq:Psi_int_exp_over_acc} and \eqref{eq:tilde_ups_aprox} hold true not only for $\cK_\infty$ but also for $\cK_\bj$.
It can be shown using an identical argument.
Therefore, we deduce from \eqref{eq:Psi_kj_bound} that
\begin{align}
    \label{eq:tilde_ups_kj_bound}
    \notag
    \left\|\breve{\Psi}_\bj \right\|_{W^{m, \infty}(\cK_\bj)}
    &\leq \left\|\Psi_\bj[\psi]\right\|_{W^{m, \infty}(\cK_\bj)}
    + \left\|\breve{\Psi}_\bj - \integral{\cU_\bj}\psi(u)\overline{\exp}(-\sR_\bj^\top a_\bj(u)) \dd u \right\|_{W^{m, \infty}(\cK_\bj)} \\
    \notag
    &\quad + \left\|\Psi_\bj[\psi] - \integral{\cU_\bj}\psi(u)\overline{\exp}(-\sR_\bj^\top a_\bj(u)) \dd u \right\|_{W^{m, \infty}(\cK_\bj)} \\
    &\leq \exp\{\cO(m \log(m D R \sigma^{-2}))\}\eps^d .
\end{align}
We also note that the derived bound remains valid if $m$ is reset to $m + 1$.
Now from  \eqref{eq:phi_j_clip_acc}, \eqref{eq:tilde_ups_kj_bound} and Lemma \ref{lem:comp_sob_norm} we find that
\begin{align*}
    &\|\breve{\Psi}_\bj - \breve{\Psi}_\bj \circ \varphi_{\bj, clip}\|_{W^{m, \infty}(\cK_\infty)} \\
    &\quad \leq \exp\{\cO(m \log(mD))\} \|\breve{\Psi}_\bj\|_{W^{m + 1, \infty}(\cK_\bj)} \eps_0 \max_{1 \leq i \leq \dim a_\bj} (1 \vee (\|(\id)_i\|_{W^{m, \infty}(\cK_\infty)} + \eps')^{2m}) \\
    &\quad \leq \exp\{\cO(m \log(m D R \sigma^{-2}))\}\eps^d \eps_0 .
\end{align*}
As a result, from \eqref{eq:Psi_int_exp_over_acc} and \eqref{eq:tilde_ups_aprox} it follows that
\begin{align}
    \label{eq:breve_Psi_phi_clip_acc}
    \|\breve{\Psi}_\bj \circ \varphi_{\bj, clip} - \Psi_\bj[\psi] \|_{W^{m, \infty}(\cK_\infty)}
    \leq \exp\{\cO(m \log(m D R \sigma^{-2}))\}\eps^d \eps_0 .
\end{align}
Meanwhile, from \eqref{eq:phi_j_clip_R}, \eqref{eq:tilde_ups_kj_bound} and Lemma \ref{lem:comp_sob_norm} we deduce that
\begin{align}
    \label{eq:breve_Psi_phi_clip_bound_R}
    \notag
    \|\breve{\Psi}_\bj \circ \varphi_{\bj, clip}\|_{W^{m, \infty}(\R^D)}
    &\leq \exp\{\cO(m \log (mD + \max_{1 \leq i \leq D} \|(\varphi_{\bj, clip})_i\|_{W^{m, \infty}(\R^D)} ))\} \| \breve{\Psi}_\bj \|_{W^{m, \infty}(\cK_\bj)} \\
    &\leq \exp\{\cO(m^2 \log(mDR \sigma^{-2}\log(1 / \eps_0)) )\} \eps^d .
\end{align}
Next, we specify the configuration of $\varphi_{\bj, clip} \in \NN(L_{clip}, W_{clip}, S_{clip}, B_{clip})$ using Lemmata \ref{lem:clip_gelu_approx} and \ref{lem:paral_nn}:
\begin{align*}
    L_{clip} \lesssim 1, \quad \|W_{clip}\|_\infty \vee S_{clip} \lesssim D,
    \quad \log B_{clip} \lesssim \log(mDR) + \log(1 / \eps_0) .
\end{align*}
Finally, we conclude from the above result, \eqref{eq:cfg_tilde_ups} and Lemma \ref{lem:concat_nn} that $\breve{\Psi}_\bj \circ \varphi_{\bj, clip} \in \NN(\breve{L}_\bj, \breve{W}_\bj, \breve{S}_\bj, \breve{B}_\bj)$ with
\begin{align}
    \label{eq:Psi_j_phi_clip_cfg}
    \notag
    &\breve{L}_\bj \lesssim \log(m + \log(1 / \eps_0)),
    \quad \|\breve{W}_\bj\|_\infty \vee \breve{S}_\bj \lesssim D^2 (m + \log(1 / \eps_0))^{3 + P(d, \beta)} \\
    &\log \breve{B}_\bj \lesssim m^5 \log^3(1 / \eps_0) \log(mDR\sigma^{-2}) ,
\end{align}

\noindent
\textbf{Step 3: increasing the number of layers.}\quad
To make the number of layers constant across all $\bj \{1, \dots, N\}^d$, we stack identity layers on top of each $\breve{\Psi}_\bj \circ \varphi_{\bj, clip}$ if necessary.
For this purpose, let $\varphi_{id, \bj}$ be the identity approximation from Lemma \ref{lem:id_deep_gelu_approx} with the accuracy parameter $\eps_0$, the scale parameter $\eps^{-d}\|\breve{\Psi}_\bj \circ \varphi_{\bj, clip}\|_{W^{m, \infty}(\R^D)}$ and the number of layers
\begin{align*}
     \log(m + \log(1 / \eps_0)) \vee \max_{\bj' \in \{1, \dots, N\}^d }\breve{L}_{\bj'} - \breve{L}_\bj + 1.
\end{align*}
We also note that the number of resulting layers of each network is chosen in such a way to achieve the upper bound for $\breve{L}_\bj$ given in \eqref{eq:Psi_j_phi_clip_cfg}.
Without loss of generality, we may assume that at least two layers are added, since otherwise the network remains unchanged.
Therefore, from \eqref{eq:breve_Psi_phi_clip_bound_R} and Lemma \ref{lem:comp_sob_norm} we find that for $\tilde{\Psi}_\bj = \eps^d\varphi_{id, \bj} \circ (\eps^{-d}\breve{\Psi}_\bj) \circ \varphi_{\bj, clip}$ the following holds:
\begin{align*}
    \|\tilde{\Psi}_\bj - \breve{\Psi}_\bj \circ \varphi_{\bj, clip} \|_{W^{m,\infty}(\cK_\infty)}
    &\leq \exp\{\cO(m \log(mD))\}\eps^d\eps_0(1 \vee \|\eps^{-d}(\breve{\Psi}_\bj \circ \varphi_{\bj, clip})\|_{W^{m, \infty}(\cK_\infty)}^m) \\
    &\leq \exp\{\cO(m^3 \log(mDR \sigma^{-2}\log(1 / \eps_0)) )\} \eps^d \eps_0 .
\end{align*}
Combining this bound and \eqref{eq:breve_Psi_phi_clip_acc}, we obtain that
\begin{align}
    \label{eq:tilde_Psi_acc_K}
    \|\tilde{\Psi}_\bj - \Psi_\bj[\psi] \|_{W^{m, \infty}(\cK_\infty)}
    \leq \exp\{\cO(m^3 \log(mDR \sigma^{-2}\log(1 / \eps_0)) )\} \eps^d \eps_0 .
\end{align}
Therefore, choosing
\begin{align*}
    \log(1 / \eps_0) \asymp m^3 \log(1 / \eps') + m^3\log(mDR\sigma^{-2})
\end{align*}
ensures that
\begin{align*}
    \|\tilde{\Psi}_\bj - \Psi_\bj[\psi] \|_{W^{m, \infty}(\cK_\infty)}
    \leq \eps^d \eps' .
\end{align*}
Furthermore, from \eqref{eq:breve_Psi_phi_clip_bound_R}, Lemma \ref{lem:id_deep_gelu_approx} and Lemma \ref{lem:comp_sob_norm} we deduce that
\begin{align}
    \label{eq:tilde_Psi_R}
    \notag
    \|\tilde{\Psi}_\bj\|_{W^{m, \infty}(\R^D)}
    &\leq \exp\{\cO(m\log (mD))\} \eps^d \|\varphi_{id, \bj}\|_{W^{m, \infty}(\R)} (1 \vee \|\eps^{-d}(\breve{\Psi}_\bj \circ \varphi_{\bj, clip})\|_{W^{m, \infty}(\R^D)}^m) \\
    &\leq \exp\{\cO(m^3 \log(mDR \sigma^{-2}\log(1 / \eps')) )\} \eps^d .
\end{align}
From \eqref{eq:Psi_j_phi_clip_cfg}, Lemmata \ref{lem:id_deep_gelu_approx} and \ref{lem:concat_nn} we conclude that $\tilde{\Psi}_\bj \in \NN(L, W, S, B)$ with
\begin{align*}
    &L \asymp \log(mDR\sigma^{-2}) + \log\log(1 / \eps'),
    \quad \log B \lesssim m^{14} \log^3(1 / \eps') \log^4(mDR\sigma^{-2}), \\
    &\|W\|_\infty \vee S \lesssim D^2 ( m^3 \log(1 / \eps') + m^3\log(mDR\sigma^{-2}) )^{3 + P(d, \beta)}, 
\end{align*}
thereby completing the proof.

\myendproof

\section{Approximation properties of GELU neural networks}

This section summarizes the approximation results presented in \cite{yakovlev2025gelu}.
Furthermore, it includes results on neural network concatenation and parallelization, alongside an evaluation of Sobolev norms for concatenation and function multiplication that are also provided in \cite{yakovlev2025gelu}.

\begin{Lem}[approximation of identity operation with multiple layers (\cite{yakovlev2025gelu}, Lemma 3.2)]
    \label{lem:id_deep_gelu_approx}
    Let $m \in \N$ and let $\id : x \mapsto x$.
    Then, for every $\eps \in (0, 1)$, every $L \in \N$ with $L \geq 2$, and every $K \geq 1$, there exists $\varphi_{id} \in \NN(L, W, S, B)$ such that
    \begin{align*}
        (i) &\quad \|\varphi_{id} - \id\|_{W^{m, \infty}([-K, K])} \leq \eps, \\
        (ii) &\quad \|\varphi_{id}\|_{W^{m, \infty}(\R)} \leq \exp\{\cO(m \log(m\log(1 / \eps)) + \log(2K))\}.
    \end{align*}
    Moreover, it holds that
    \begin{align*}
        \|W\|_\infty \lesssim 1, \quad S \lesssim L,
        \quad \log B \lesssim (m + L)\log m + \log(1 / \eps) + m\log(K).
    \end{align*}
\end{Lem}

\begin{Lem}[approximation of clipping operation (\cite{yakovlev2025gelu}, Lemma 3.5)]
    \label{lem:clip_gelu_approx}
    For every $A \geq 1$, every $\eps \in (0, 1)$, and every $m \in \N$, there exists $\varphi_{clip} \in \NN(L, W, S, B)$ such that
    \begin{align*}
        (i) &\quad \|\varphi_{clip} - \id \|_{W^{m, \infty}([-A, A])} \leq \eps, \\
        (ii) &\quad \|\varphi_{clip} + A + 1/2\|_{W^{m, \infty}((-\infty, -A - 1])} \vee \|\varphi_{clip} - A - 1/2\|_{W^{m, \infty}([A + 1, +\infty))}
        \leq \eps, \\
        (iii) &\quad \|\varphi_{clip}\|_{W^{m, \infty}(\R)} \leq \exp\{\cO(m \log m + m\log\log(1 / \eps) + \log(2A) )\}, \\
        (iv) &\quad \|\varphi_{clip}\|_{W^{0, \infty}(\R)} \leq A + 5 / 2, \\
        (v) &\quad \|\varphi_{clip} + A + 1/2\|_{W^{0, \infty}((-\infty, -A])} \vee \|\varphi_{clip} - A - 1/2\|_{W^{0, \infty}([A, +\infty))} \leq \eps + 1, \\
        (vi) &\quad |\varphi_{clip}|_{W^{k, \infty}(\R)} \leq \exp\{\cO(k \log m + k \log\log(1 / \eps) )\} .
    \end{align*}
    Moreover, $\varphi_{clip}$ has $L = \|W\|_\infty = 2$, $S = 7$ and $\log B \lesssim \log(A m / \eps)$.
    
\end{Lem}

\begin{Lem}[approximation of square operation (\cite{yakovlev2025gelu}, Lemma 3.6)]
\label{lem:square_approx}
    Define $f_{sq} : x \mapsto x^2$.
    Then, for every $\eps \in (0, 1)$ and every $m \in \N$, there exists $\varphi_{sq} \in \NN(L, W, S, B)$ such that
    \begin{align*}
        \|\varphi_{sq} - f_{sq}\|_{W^{m, \infty}([-C, C])} \leq C^3 \eps, \quad \text{for all } C \geq 1 .
    \end{align*}
    Furthermore, $L = \|W\|_\infty = 2$, $S = 6$ and $\log B \lesssim  \log(1 / \eps) + \log m$.
\end{Lem}

\begin{Lem}[approximation of two number multiplication (\cite{yakovlev2025gelu}, Corollary 3.7)]
    \label{lem:multi_approx_gelu}
    Define $\mathrm{prod}_2 : (x, y) \mapsto x \cdot y$, and let $m \in \N$ be arbitrary.
    Then, for any $\eps \in (0, 1)$, there exists $\varphi_{mul} \in \NN(L, W, S, B)$ satisfying
    \begin{align*}
        \|\varphi_{mul} - \mathrm{prod}_2\|_{W^{m, \infty}([-C, C]^2)}
        \leq C^3 \eps ,
        \quad \text{for all } C \geq 1 .
    \end{align*}
    In addition, $L = 2$, $\|W\|_\infty \leq 4$, $S \leq 12$ and $\log B \lesssim \log (1/\eps) + \log m$.
\end{Lem}

\begin{Lem}[approximation of multivariate polynomials (\cite{yakovlev2025gelu}, Lemma 3.10)]
    \label{lem:part_sum_monomial}
    Define $f_{\mathcal{A}} : x \mapsto \sum_{\bk \in \mathcal{A} } a_\bk x^\bk$, where $x \in \R^I$ and $\mathcal{A} = \{\bk \in \Z_+^I \; : \; |\bk| \leq d \}$ for some $I, d \in \N$ with $d \geq 2$.
    Also assume that $|a_\bk| \leq 1$ for all $\bk \in \mathcal{A}$.
    Then, for every $\eps \in (0, 1)$, every natural $m \geq 3$ and every $K \geq 1$, there exists a neural network $\varphi_{\mathcal{A}} \in \NN(L, W, S, B)$ such that
    \begin{align*}
        (i) &\quad \|f_\mathcal{A} - \varphi_\mathcal{A}\|_{W^{m, \infty}([-K, K]^I)} \leq \eps, \\
        (ii) &\quad \|\varphi_\mathcal{A}\|_{W^{m, \infty}(\R^I)} \leq \exp\{\cO( (m^2 + md + I)\log(mdKI\log(1 / \eps)) )\}.
    \end{align*}
    In addition, $\varphi_{\mathcal{A}}$ has
    \begin{align*}
        &L \lesssim \log d,
        \quad \|W\|_\infty \vee S \lesssim (d + I)^{3 + d \wedge I} \\
        &\log B \lesssim (\log(1 / \eps) + m^2(d + I)\log(mdKI) + m^2d^2)\log(d + I).
    \end{align*}    
\end{Lem}

\begin{Lem}[approximation of the exponential function (\cite{yakovlev2025gelu}, Lemma 3.11)]
    \label{lem:min_exp_gelu_approx}
    Define $f_{exp} : x \mapsto e^{-x}$, and let $m \in \N$ be arbitrary.
    Then, for any $\eps \in (0, 1)$ and $0 \leq A \leq 1$, there exists a neural network $\varphi_{exp} \in \NN(L, W, S, B)$ such that
    \begin{align*}
        (i) &\quad \|\varphi_{exp} - f_{exp}\|_{W^{m, \infty}([-A, +\infty))} \leq \eps, \\
        (ii) &\quad \|\varphi_{exp}\|_{W^{m, \infty}(\R)} \leq \exp\{\cO(m^2\log(m\log(1 / \eps))  )\} .
    \end{align*}
    Furthermore, 
    \begin{align*}
        L \lesssim \log m + \log\log(1 / \eps),
        \quad \|W\|_\infty \vee S \lesssim m^{12} \log^4(1 / \eps), 
        \quad \log B \lesssim m^{11} \log^3(1 / \eps).
    \end{align*}
\end{Lem}

\begin{Lem}[partition of unity approximation (\cite{yakovlev2025gelu}, Lemma 3.4)]
    \label{lem:pou_gelu_approx}
    Define $a_i = 2^{-N + i}$ for each $i \in \{0, 1, \dots, N\}$, where $N \in \N$ and $N \geq 3$.
    Then, for every $\eps \in (0, 1)$ and every $m \in \N$, there exist $\{\psi_i\}_{i=1}^N$, with $\psi_i \in \NN(L, W, S, B)$ for each $1 \leq i \leq N$, such that
    \begin{align*}
        (i) &\quad \sum_{i=1}^N \psi_i(x) = 1, \quad \text{for all } x \in \R, \\
        (ii) &\quad \max_{1 \leq i \leq N} \|\psi_i\|_{W^{m, \infty}(\R)} \leq \exp\{\cO(mN + m\log(m\log(1 / \eps)))\}, \\
        (iii) &\quad \|\psi_N\|_{W^{m, \infty}(-\infty, a_{N - 2}]} \vee \|\psi_1\|_{W^{m, \infty}([a_2, +\infty))} \vee  \max_{2 \leq i \leq N - 1}\|\psi_i\|_{W^{m, \infty}(\R \setminus (a_{i - 2}, a_{i + 1}))}
        \leq \eps,
    \end{align*}
    Furthermore, $L = 2$, $\|W\|_\infty \vee S \lesssim 1$ and $\log B \lesssim \log(1 / \eps) + mN + m\log m$.
\end{Lem}

\begin{Lem}[division operation approximation (\cite{yakovlev2025gelu}, Lemma 3.14)]
    \label{lem:div_gelu_approx}
    Define $\mathrm{div} : (x, y) \mapsto x / y$ for any $x \in \R$ and $y > 0$.
    Let also $a_0 = 2^{-N}$ for $N \in \N$ with $N \geq 3$.
    Then, for every $\eps \in (0, 1)$ and every $m \in \N$ such that $m \geq 3$, there exists a GELU network $\varphi_{div} \in \NN(L, W, S, B)$ satisfying
    \begin{align*}
        (i) &\quad \|\varphi_{div} - \mathrm{div}\|_{W^{m, \infty}([-1, 1] \times [a_0, 1])} \leq \eps, \\
        (ii) &\quad \| \varphi_{div} \|_{W^{m, \infty}(\R^2)} \leq \exp\{\cO(m^4 N + m^4\log(m\log(1 / \eps)))\} .
    \end{align*}
    Furthermore, the network $\varphi_{div}$ has
    \begin{align*}
        L \lesssim \log(mN\log(1 / \eps)),
        \quad \|W\|_\infty \vee S \lesssim m^{21} N^5 \log^4(1 / \eps),
        \quad \log B \lesssim m^{24} N^4 \log^4(1 / \eps).
    \end{align*}    
\end{Lem}

\begin{Lem}[properties of GELU activation function (\cite{yakovlev2025generalization}, Lemma B.2)]
\label{lem:gelu_seminorms_bound}
    For any $k \in \N$ we have the following bounds for the Sobolev seminorms:
    \begin{align*}
        \left|\gelu\right|_{W^{k, \infty}(\R)} \leq
        \begin{cases}
            1 + 1 / \sqrt{2\pi}, \quad &k = 1, \\
            (k + 1)\sqrt{\frac{(k - 2)!}{2\pi}}, \quad &k \geq 2
        \end{cases}
    \end{align*}
    For $k = 0$ we have that
    \begin{align*}
        \|\gelu\|_{W^{0, \infty}([-C, C])} \leq C, \quad \text{for all } C > 0.
    \end{align*}
    In addition, for any $A \geq 0$ and $m \in \N$, the tails behave as follows:
    \begin{align*}
        \|\gelu - \id\|_{W^{m, \infty}([A, +\infty))} \vee \|\gelu\|_{W^{m, \infty}((-\infty, -A])}
        \leq 2 e^{-A^2 / 4} \sqrt{m!} .
    \end{align*}
\end{Lem}

\begin{Lem}[\cite{de2021approximation}, Lemma A.6]
    \label{lem:prod_sob_norm}
    Let $d \in \N$, $k \in \Z_+$, $\Omega \subseteq \R^d$ and $f, g \in W^{k, \infty}(\Omega)$.
    Then it holds that
    \begin{align*}
        \|f \cdot g\|_{W^{k, \infty}(\Omega)} \leq 2^k \|f\|_{W^{k, \infty}(\Omega)} \|g\|_{W^{k, \infty}(\Omega)}.
    \end{align*}
\end{Lem}

\begin{Lem}[\cite{yakovlev2025gelu}, Lemma B.4]
    \label{lem:comp_sob_norm}
    Let $d, m, n \in \N$ and let also $\Omega_1 \subseteq \R^d$, $\Omega_2 \subseteq \R^m$, $f \in C^n(\Omega_1, \Omega_2)$ and $g \in C^n(\Omega_2, \R)$. Then it holds that
    \begin{align*}
        \|g \circ f\|_{W^{n, \infty}(\Omega_1)} \leq 16 (e^2n^4md^2)^n \|g\|_{W^{n, \infty}(\Omega_2)} \max_{1 \leq i \leq m}(\|(f)_i\|_{W^{n, \infty}(\Omega_1)}^n \vee 1).
    \end{align*}
    Moreover, if $g \in C^{n + 1}(\Omega_2, \R)$ and $\tilde f \in C^n(\Omega_1, \Omega_2)$, then
    \begin{align*}
        &\|g \circ f - g \circ \tilde{f}\|_{W^{n, \infty}(\Omega_1)}
        \\&
        \leq 32(e^2n^5 m^2 d^2)^n \|g\|_{W^{n + 1, \infty}(\Omega_2)} \max_{1 \leq i \leq m}\|(f)_i - (\tilde{f})_i\|_{W^{n, \infty}(\Omega_1)} \left(1 \vee \|(f)_i\|^{2n}_{W^{n, \infty}(\Omega_1)} \vee \|(\tilde{f})_i\|^{2n}_{W^{n, \infty}(\Omega_1)}\right).
    \end{align*}
\end{Lem}

\begin{Lem}[concatenation of neural networks (\cite{yakovlev2025gelu}, Lemma B.5)]
    \label{lem:concat_nn}
    Given $K \in \N$ with $K \geq 2$. Then, for any neural networks $\varphi^{(k)} \in \NN(L_k, W_k, S_k, B_k)$ with $1 \leq k \leq K$
    such that $\varphi^{(k)} : \R^{d_k} \to \R^{d_{k + 1}}$, there exists a neural network $h = \varphi^{(K)} \circ \varphi^{(K - 1)} \dots \circ \varphi^{(1)} \in \NN(L, W, S, B)$ satisfying
    \begin{align*}
        L &\leq 1 + \sum_{k=1}^K (L_k - 1), \quad S \leq \sum_{k=1}^K S_k + 2\sum_{k = 1}^{K - 1}\|W_k\|_\infty \cdot \|W_{k + 1}\|_\infty, \\
        \|W\|_{\infty} &\leq \max_{1 \leq k \leq K}\|W_k\|_{\infty},
        \quad B \leq 2 \max_{1 \leq k \leq K - 1} \left[ (B_k \vee 1) (B_{k + 1} \vee 1) \left( \|W_k\|_\infty \vee \|W_{k + 1}\|_\infty \right) \right].
    \end{align*}    
\end{Lem}

\begin{Lem}[parallelization of neural networks (\cite{yakovlev2025gelu}, Lemma B.6)]
    \label{lem:paral_nn}
    Let $K \in \N$ with $K \geq 2$ and let neural networks $\varphi^{(k)} \in \NN(L_k, W_k, S_k, B_k)$ for $1 \leq k \leq K$.
    Assume further that $L_k = L$ for all $1 \leq k \leq K$.
    Then, the following holds:
    \begin{itemize}
        \item[(i)] if $\varphi^{(k)} : \R^{d_k} \to \R$ for each $1 \leq k \leq K$, then there exists a neural network $\varphi \in \NN(L, W, S, B)$ such that $\varphi(x) = (\varphi^{(1)}(x_1), \dots, \varphi^{(K)}(x_K))^\top$ for all $x = (x_1^\top, \dots, x_K^\top)^\top$,
        where $x_k \in \R^{d_k}$ for any $1 \leq k \leq K$.
        In addition, there exists $\varphi_{sum} \in \NN(L, W, S, B_{sum})$, which implements the summation, that is,
        \begin{align*}
            \varphi_{sum}(x) = \sum_{k=1}^K \varphi^{(k)}(x_k), \quad \text{for all } x = (x_1^\top, \dots, x_K^\top)^\top .
        \end{align*}
        \item[(ii)] if $\varphi^{(k)} : \R^p \to \R$ for some $p \in N$ for every $1 \leq k \leq K$, then there exists a neural network $\varphi \in \NN(L, W, S, B)$ satisfying $\varphi(x) = (\varphi^{(1)}(x), \dots, \varphi^{(K)}(x))^\top$ for all $x \in \R^p$.
        Moreover, there exists a summation network $\varphi_{sum} \in \NN(L, W, S, B_{sum})$ such that
        \begin{align*}
            \varphi_{sum}(x) = \sum_{k=1}^K \varphi^{(k)}(x), \quad \text{for all } x \in \R^p .
        \end{align*}
    \end{itemize}
    Furthermore, in both cases it holds that
    \begin{align*}
        \|W\|_\infty \leq \sum_{k=1}^K \|W^{(k)}\|_\infty,
        \quad S \leq \sum_{k = 1}^K S^{(k)},
        \quad B \leq \max_{1 \leq k \leq K} B^{(k)},
        \quad B_{sum} \leq K \max_{1 \leq k \leq K} B^{(k)} .
    \end{align*}
\end{Lem}

\section{Auxiliary results}

This section provides auxiliary approximation results from \cite{yakovlev2025generalization} and relevant tools from probability theory presented in \cite{wainwright19}.

\begin{Lem}[\cite{yakovlev2025generalization}, Lemma A.1]
    \label{lem:g_circ_L_inf_acc}
    Let $g^* \in \cH^\beta([0, 1]^d, \R^D, H)$ and let $g^\circ$ be as defined in \eqref{eq:g_circ_def}.
    Then it holds that
    \begin{align*}
        \|g^* - g^\circ\|_{L^\infty([0, 1]^d)}
        \leq \frac{H d^\floor{\beta}\eps^\beta\sqrt{D}}{\floor{\beta}!} .
    \end{align*}
    
\end{Lem}

\begin{Lem}[\cite{wainwright19}, Proposition 2.9]
    Suppose that a random variable $X$ is sub-exponential with parameters $(\nu, \alpha)$, that is
    \begin{align*}
        \E \exp\{\lambda(X - \E X)\} \leq \exp\{\nu^2\lambda^2 / 2\}, \quad \text{for all} |\lambda| \leq 1 / \alpha .
    \end{align*}
    Then, for all $t \geq 0$, it holds that
    \begin{align*}
        \p(X \geq \E[X] + t) \leq \exp\left\{-\frac{1}{2}\left(\frac{t^2}{\nu^2} \wedge \frac{t}{\alpha}\right)\right\}
    \end{align*}

\end{Lem}

\begin{Lem}[\cite{wainwright19}, Example 2.11]
    \label{lem:chi_sq_concentr}
    A chi-squared random variable with $D$ degrees of freedom is sub-exponential with parameters $(\nu, \alpha) = (2\sqrt{D}, 4)$.
    In particular, for $X \sim \chi^2(D)$ and for all $t \geq 0$, it holds that
    \begin{align*}
        \p(X \geq D + t) \leq \exp\left\{-\frac{1}{8}\left(\frac{t^2}{D} \wedge t \right)\right\} .
    \end{align*}
\end{Lem}

\end{document}